\documentclass[12pt]{amsart}

\usepackage{tikz}

\usepackage[latin1]{inputenc}
\usepackage{amsmath} 
\usepackage{amsfonts}
\usepackage{amssymb}
\usepackage{latexsym}
\usepackage{graphicx}
\usepackage{subfigure}
\usepackage{color}
\usepackage{hyperref}
\usepackage{verbatim}
\usepackage[all]{xy}
\usepackage{graphics}
\usepackage{pdfsync}

\DeclareMathOperator{\tr}{tr}
\DeclareMathOperator{\Tb}{Tb}
\DeclareMathOperator{\s}{s}

\theoremstyle{plain}
\newtheorem{theorem}{Theorem}[section]
\newtheorem{lemma}[theorem]{Lemma}
\newtheorem{corollary}[theorem]{Corollary}

\newtheorem{definition}[theorem]{Definition}
\newtheorem{example}[theorem]{Example}
\newtheorem{proposition}[theorem]{Proposition}

\theoremstyle{remark}
\newtheorem{remark}[theorem]{Remark}

\title[Periodic P-Partitions]{Periodic P-Partitions}

\author{Brian T. Chan}
\email{bchan600@gmail.com}

\date{\today}


\newlength\cellsize 
\savebox2{%
\begin{picture}(15,15)
\put(0,0){\line(1,0){15}}
\put(0,0){\line(0,1){15}}
\put(15,0){\line(0,1){15}}
\put(0,15){\line(1,0){15}}
\end{picture}}
\newcommand\cellify[1]{\def\thearg{#1}\def\nothing{}%
\ifx\thearg\nothing
\vrule width0pt height\cellsize depth0pt\else
\hbox to 0pt{\usebox2\hss}\fi%
\vbox to 15\unitlength{
\vss
\hbox to 15\unitlength{\hss$#1$\hss}
\vss}}
\newcommand\tableau[1]{\vtop{\let\\=\cr
\setlength\baselineskip{-16000pt}
\setlength\lineskiplimit{16000pt}
\setlength\lineskip{0pt}
\halign{&\cellify{##}\cr#1\crcr}}}
\savebox3{%
\begin{picture}(15,15)
\put(0,0){\line(1,0){15}}
\put(0,0){\line(0,1){15}}
\put(15,0){\line(0,1){15}}
\put(0,15){\line(1,0){15}}
\end{picture}}
\newcommand\expath[1]{%
\hbox to 0pt{\usebox3\hss}%
\vbox to 15\unitlength{
\vss
\hbox to 15\unitlength{\hss$#1$\hss}
\vss}}
\newcommand\bas[1]{\omit \vbox to \cellsize{ \vss \hbox to \cellsize{\hss$#1$\hss} \vss}}

\begin{document}

\begin{abstract} In this paper, we introduce a class of $(P, \omega)$-partitions that we call \emph{periodic $(P, \omega)$-partitions}, then prove that such $(P, \omega)$-partitions satisfy a homogeneous first-order matrix difference equation. After defining an appropriate counting problem for the above $(P, \omega)$-partitions, we show that as a consequence of this equation, periodic $(P, \omega)$-partitions can be enumerated with constant coefficient linear recurrence relations. By analysing the above matrix difference equation, we also prove a result for the asymptotic growth rate for the number of periodic $(P, \omega)$-partitions. The results of this paper generalizes and strengthens the constant coefficient linear recurrence results proved by Sun and by L\'opez, Mart\'inez, P\'erez, P\'erez, and Basova for enumerating standard Young tableaux on shifted strips with constant width.
\end{abstract}
\maketitle
%
\section{Introduction} \label{sec:intro}

The problem of enumerating standard Young tableaux on non-classical shapes is a topic of recent interest \cite{SYT} in which shapes other than a Young diagram or a skew shape are considered. For instance, a problem equivalent to the problem of enumerating the number of standard Young tableaux of shifted strips was considered by Hardin and Heinz (\cite{OEIS}, A181196). Hardin and Heinz (\cite{OEIS}, A181196; cf. \cite{OEIS}, A227578) conjectured constant coefficient linear recurrence relations for enumerating the number of standard Young tableaux on shifted strips with constant width up to when the width is seven. Sun \cite{EYTPS} proved Hardin and Heinz's relations for when the constant width is four and five. Moreover, L\'opez, L\'opez, Mart\'inez, P\'erez, P\'erez, and Basova \cite{CHC}, calling shifted strips parallelogramic shapes, generalized Sun's constant coefficient linear recurrence results to all shifted strips with a fixed number of cells in each row. Tewari and van~Willigenburg proved that the above enumeration results are connected to the representation theory of the $0$-Hecke algebra \cite{HAlgebra} and proved Hardin and Heinz' conjecture in the case where the constant width is three. \\

$(P, \omega)$-partitions were first considered by MacMahon in \cite{CA}. Later on, the theory of $(P, \omega)$-partitions was developed by Gessel and Stanley \cite{MPPIP, StanThesis}. For instance, this theory is known to have applications to quasisymmetric functions as $(P, \omega)$-partitions are essential for the theory of quasisymmetric functions \cite{MPPIP}. In this paper, we investigate a different aspect of $(P, \omega)$-partitions by generalizing the results described in the previous paragraph to certain $(P, \omega)$-partitions, then prove an asymptotic property of the sequences we are enumerating and establish a divisibility property relating $(P, \omega)$-partitions to minimal polynomials of certain matrices. \\

Specifically, we introduce a class of $(P, \omega)$-partitions that we call periodic $(P, \omega)$-partitions, define an appropriate counting problem in order to enumerate such $(P, \omega)$-partitions, and generalize the coefficient linear recurrence results described in the above paragraph to this class. To accomplish this, we introduce the notion of a connected triple for a poset to prove the existence of a matrix difference equation that can enumerate the above class of $(P, \omega)$-partitions. As consequences of our enumeration result, we give new proofs of Sun's and L\'opez et.al.'s constant coefficient linear recurrence results and explain how semi-standard variants of Sun and L\'opez's results can also be derived. We also give a brief outline of how our results can also be applied to standard Young tableaux on certain truncated shifted shapes that exhibit a fixed repeating pattern and for higher-dimensional analogues of such tableaux. \\

We then prove that the matrix in the matrix difference equation we construct, that we call a tableau transfer matrix, is a primitive matrix in many cases. Using this result, we apply Perron-Frobenius theory to the tableau transfer matrices that are primitive to prove a result for the asymptotic growth rate of the number of periodic $(P, \omega)$-partitions. Lastly, we briefly remark on how our proof that these matrices are primitive can be used to give further descriptions for the constant coefficient linear recurrence relations we are analysing. \\

In Section \ref{sec:prelims}, we give the preliminaries and describe the appropriate counting problem that we will use to enumerate the $(P, \omega)$-partitions. In Section \ref{sec:periodic}, we introduce connected triples and periodic $(P, \omega)$-partitions, then describe examples of such $(P, \omega)$-partitions. In Section \ref{sec:matrix}, we introduce the tableau transfer matrices and prove that the number of periodoic $(P, \omega)$-partitions satisfy a matrix difference equation, then we prove that the number of periodic $(P, \omega)$-partitions can be enumerated by a constant coefficient linear recurrence relation and, as corollaries, give new proofs of Sun and L\'opez et.al's results. Lastly, in Section \ref{sec: asymptotic}, we prove a result for the asymptotic growth rate of the number of periodic $(P, \omega)$-partitions.

\section{Preliminaries and the Enumeration Problem} \label{sec:prelims}

Let $\mathbb{N}$ denote the set of positive integers and let $\mathbb{N}_0$ denote the set of non-negative integers. If $f : X \to Y$ is a function and if $X'$ is a subset of $X$, then let $f|_{X'}$ denote the restriction of $f$ to $X'$. Let $X$ be a set, and let $\equiv$ be an equivalence relation on $X$. Then for all $r \in X$, the \emph{equivalence class of $r$ in $X$ with respect to $\equiv$} is the set $\{r_1 \in X : r \equiv r_1\}$. Lastly, let $X / \equiv$ denote the set of all equivalence classes in $X$ with respect to $\equiv$. If $X$ is a (non-empty) set, then a \emph{set partition of $X$} is a set $\mathcal{F}$ of non-empty subsets of $X$ such that every element of $X$ is contained in exactly one element of $\mathcal{F}$. For all positive integers $n$, define $[n] = \{1,2,\dots,n\}$. Moreover, for all positive integers $n_1$ and $n_2$ such that $n_1 \leq n_2$, define $[n_1,n_2] = \{k \in \mathbb{N} : n_1 \leq k \leq n_2 \}$. In particular, $[n,n] = \{n\}$, $[n,n+1] = \{n, n+1\}$, $[n,n+2] = \{n, n+1, n+2\}$, and so on. Lastly, write $[n, \infty) = \{m \in \mathbb{Z} : m \geq n\}$. We denote any sequence of numbers $a_1$, $a_2$, $\dots$ by $(a_n)_{n=1,2,\dots}$. Such a sequence satisfies a \emph{constant coefficient linear recurrence relation} if there exist constants $c_1$, $c_2$, $\dots$, $c_d$ such that for all $n \geq 1$, $a_{n+d} = c_1 a_{n + d - 1} + c_2 a_{n + d - 2} + \cdots + c_d a_n$. The \emph{characteristic polynomial} of the aforementioned recurrence $a_{n+d} = c_1 a_{n + d - 1} + c_2 a_{n + d - 2} + \cdots + c_d a_n$ is $x^d - c_1 x^{d-1} - c_2 x^{d-2} - \cdots - c_d$. \\

We will use posets extensively in this paper, so we adhere to the following conventions. Recall the notion of \emph{order-preserving maps}, \emph{order-reversing maps}, \emph{order-isomorphisms}, \emph{order-embeddings}, and \emph{subposets} \cite{ILO, OS}. Moreover, call an \emph{order-automorphism on $P$} \cite{OS} an order-isomorphism of the form $f : P \to P$ . For convenience, we write $P$ instead of $(P, \leq)$. Whether $P$ represents a poset or its underlying set of elements will be clear from context. For instance, an order-preserving map will be written as $f : P \to Q$. Subsets of the set of elements of a poset $P$ will be regarded as subposets of $P$ and vice versa. In particular, given a poset $P$, we may also call a subposet of $P$ a subset of $P$. We regard any poset as being a subposet of itself, and we allow for a poset to have no elements. \\

As we regard subposets as subsets and vice versa, we will also use set-theoretic notation for subposets. Given a poset $P$ and two subposets $S_1$ and $S_2$ of $P$, let $S_1 \cap S_2$ denote the subposet of $P$ whose elements are the elements of $P$ that are contained in $S_1$ and in $S_2$, let $S_1 \cup S_2$ denote the subposet of $P$ whose elements are the elements of $P$ that are contained in $S_1$ or in $S_2$, let $S_1 \backslash S_2$ denote the subposet of $P$ whose elements are the elements of $S_1$ that are not in $S_2$, and let $S_1 \subseteq S_2$ mean that every element of $S_1$ is contained in $S_2$. Moreover, if $f : P \to Q$ is a map between two posets, if $S$ is a subset of $P$, and if $S'$ is a subset of $Q$, then let $f(S)$ be the subposet of $Q$ with $\{f(p) : p \in P\}$ as its set of elements and let $f^{-1}(S')$ be the subset of $P$ with $\{p \in P : f(p) \in S'\}$ as its set of elements. Lastly, let $\emptyset$ denote the empty set and the poset with no elements. \\

Additionally, we define the following \cite{ILO}. Write $p \nleq q$ if $p \leq q$ is false. Also, let $\geq$, $>$, and $<$ take on their usual meanings. If $p, q \in P$, then write $p \parallel q$ to mean that $p \nleq q$ and $q \nleq p$. Furthermore, if we want to clarify which partial order we are using, we say $p \leq q$ in $P$, where $P$ is a poset and $p,q \in P$. Let $P$ be a poset. Then we call a subset $S$ of $P$ an \emph{order ideal} of $P$ if, for all $p \in S$ and $q \in P$, $q < p$ in $P$ implies $q \in S$. Moreover, we call a subset $S$ of $P$ an \emph{order filter} of $P$ if, for all $p \in S$ and $q \in P$, $q > p$ in $P$ implies $q \in S$. \\

A \emph{Young diagram} (cf. \cite{SYT, EC1, EC2}) is the empty set or a finite subset $X$ of $\mathbb{N}^2$ such that for some $i,j \in \mathbb{N}$, $(1,j) \in X$ and $(i,1) \in X$. We call the elements of a Young diagram $X$ the \emph{cells of $X$}. Lastly, given a Young diagram $X$, define a \emph{row of $X$} to be a non-empty subset of cells of the form $\{r \in X : \exists j \in \mathbb{N} \; \text{ such that } \; r = (i,j)\}$ for some $i \in \mathbb{N}$, and define a \emph{column of $X$} to be a non-empty subset of cells of the form $\{r \in X : \exists i \in \mathbb{N} \; \text{ such that } \; r = (i,j) \}$ for some $j \in \mathbb{N}$. Following convention, we will depict the cells of a Young diagram with boxes. Moreover, when depicting a Young diagram $X$, we position $(i,j+1)$ immediately to the right of $(i,j)$ and we position $(i+1, j)$ immediately below $(i,j)$.

\begin{example} If $X_1 = \{(1,1), (1,2), (1,3)\} $, $X_2 = \{(1,1), (1,3), (1,6), (1,7), (2,3)\}$, and $X_3 = \{(1,3), (2,2), (3,1), (3,3), (3,4) \}$, then the Young diagram $X_1$ is depicted by
\begin{center}
$\tableau{{} & {} & {}}$ \;,
\end{center}
the Young diagram $X_2$ is depicted by
\begin{center}
$\tableau{{} && {} &&& {} & {} \\ && {}}$ \;,
\end{center}
and the Young diagram $X_3$ is depicted by
\begin{center}
$\tableau{ && {} \\ & {} \\ {} && {} & {} }$ \;.
\end{center}
\end{example}

Regard $\mathbb{Z}$ as the poset $\cdots < -2 < -1 < 0 < 1 < 2 < \cdots$. Moreover, for any integer $d \geq 2$, regard the $d$-dimensional integer lattice $\mathbb{Z}^d$ as a poset with the standard partial order given by $(k_1, k_2, \dots, k_d) \leq (k'_1, k'_2, \dots, k'_d)$ if $k_i \leq k_i'$ for all $1 \leq i \leq d$. We will also regard any subset $S$ of $\mathbb{Z}^d$ as the subposet of $\mathbb{Z}^d$ with $S$ as its set of elements. In particular, as $X \subset \mathbb{Z}^2$ for any Young diagram $X$, we will also regard Young diagrams as subposets of $\mathbb{Z}^2$. In a Hasse diagram, $p < q$ implies that $q$ is positioned above $p$, but in a Young diagram, $p < q$ implies that $q$ cannot be positioned to the left of $p$ and $q$ cannot be positioned above $p$. We illustrate this with an example.

\begin{example} If $P$ is the poset depicted by the following Hasse diagram
\begin{center}
	\begin{tikzpicture}[scale=0.7]
	\node (aa) at (-2,-1) {$\circ$};
	\node (a) at (-1,0) {$\circ$};
	\node (b) at (0,-1) {$\circ$};
	\node (c) at (-1,-2) {$\circ$};
	\node (d) at (1,-2) {$\circ$};
	\node (e) at (2,-3) {$\circ$};
	\node (f) at (0,-3) {$\circ$};
	
	\draw (aa) -- (c);
	\draw (a) -- (b) -- (c) -- (b) -- (d) -- (e);
	\draw (c) -- (f) -- (d);
	\end{tikzpicture},
\end{center}
then any of the Young diagrams shown below can also be used to depict $P$.
$$\tableau{& {} & {} & {} \\ {} & {} & {} \\ && {}} \;\;\;\;\;\;\;\;\;\;\;\;\;\;\;\;\;\;\;\; \tableau{& {} \\ {} & {} \\ {} & {} & {} \\ {}}$$
\end{example}

Some posets that are not Young diagrams can also be depicted with Young diagrams. For instance, note the following. 

\begin{example} If $P$ is the subposet of $\mathbb{Z}^2$ where $P = \{(-1,-5), (-1,-4), (-1,-3)\}$, then we formally define the following Young diagram that depicts $P$ to be $\{(1,1), (1,2), (1,3)\}$.
\begin{center}
$\tableau{{} & {} & {}}$
\end{center}
\end{example}

If $u, v \in \mathbb{R}^d$ and $c_1, c_2 \in \mathbb{R}$, then let $c_1 u + c_2 v$ and $c_1 u - c_2 v$ take on their usual meanings from linear algebra where $u$ and $v$ are vectors and where $c_1$ and $c_2$ are scalars. Moreover, if $u \in \mathbb{R}^d$, if $c \in \mathbb{R}$, and if $S \subseteq \mathbb{R}^d$, then write $S + cu = \{v + cu : v \in S \}$ and write $S - cu = \{v - cu : v \in S \}$. \\

Call a matrix an \emph{integer matrix} if all entries of that matrix are integers. If $M$ is a matrix, then let $M(i,j)$ denote the entry in the $i^{th}$ row and $j^{th}$ column of $M$. A \emph{column vector} is a matrix with one column, and a \emph{row rector} is a matrix with one row. If $v$ is a column vector with $N$ rows, then for all $1 \leq i \leq N$, let $v(i)$ denote the entry in the $i^{th}$ row of $v$. If $p(x) = \sum_{k = 0}^n a_k x^k$ is a polynomial for some scalars $a_0$, $a_1$, $\dots$, $a_k$ and if $M$ is an $N$ by $N$ square matrix, then write $p(M) = a_0 I_N + \sum_{k = 1}^n a_k M^k$, where $I_N$ is the $N$ by $N$ identity matrix. Given functions $f, g : \mathbb{N} \to \mathbb{R}_{>0}$, $\mathbb{R}_{>0}$ denoting the set of positive reals, write $f(n) \sim g(n)$ as $n \to \infty$ if $\lim_{n \to \infty} f(n) / g(n) = 1$. Lastly, given sequences $(v_n)_{n=1,2,\dots}$ and $(v_n')_{n=1,2,\dots}$ of column-vectors with positive entries such that for all $n$, $v_n$ and $v_n'$ have $N$ entries, write $v_n \sim v_n'$ as $n \to \infty$ if for all integers $1 \leq i \leq N$, $v_n(i) \sim v'_n(i)$ as $n \to \infty$. \\

If $X$ is a finite set, then an \emph{indexing of $X$} is a bijection $f : X \to [k]$ where $k = |X|$. Next, we introduce terminology for $(P, \omega)$-partitions from \cite{EC1}, but extend labellings to countably infinite posets and define $\mathcal{A}(Q, \omega)$ in a non-standard way. If $P$ is a finite poset, then a \emph{labelling of $P$} is a bijection $\omega : P \to [k]$ where $k = |P|$. Moreover, if $P$ is a countably infinite poset, then a labelling of $P$ is a bijection $\omega : P \to \mathbb{Z}$. Recall \cite{EC1, EC2} that a map $f : P \to Q$, where $P$ and $Q$ are posets, is \emph{order-reversing} if $p_1 \leq p_2$ in $P$ implies $f(p_2) \leq f(p_1)$ in $Q$. If $P$ is a finite poset and if $\omega$ is a labelling of $P$, then a \emph{$(P, \omega)$-partition} is an order-reversing map $f : P \to \mathbb{N}_0$ such that $f(x) > f(y)$ if $x < y$ and $\omega(x) > \omega(y)$. We sometimes call a $(P, \omega)$-partition a \emph{$P$-partition}. Next, if $P$ is a countable poset, if $Q$ is a finite poset such that $Q \subseteq P$, and if $\omega$ is a labelling of $P$, then let $\mathcal{A}(Q, \omega)$ denote the set of order-reversing maps $U : Q \to \mathbb{N}_0$ such that $U(x) > U(y)$ if $x < y$ and $\omega(x) > \omega(y)$. If $Q_1 \subseteq Q$, then we let $U|_{Q_1}$ denote the restriction of the function $U$ to $Q_1$. We also denote labellings with the symbol $\omega'$ or with the symbol $\omega''$. Lastly, if $P$ is a countable poset and if $\omega$ is a labelling of $P$, then $\omega$ is a \emph{dual-natural labelling} if $\omega$ is order-reversing and a \emph{natural labelling} if $\omega$ is order-preserving (cf. \cite{EC1, EC2}). \\

If $P$ is a poset, then we will depict a labelling or an order-reversing map on $P$ with a Hasse diagram (or a Young diagram) whose nodes (or cells) are filled with integers. Specifically, if $P$ is a poset, if $p \in P$, if $X \subseteq \mathbb{Z}$, and if $f : P \to X$ satisfies $f(p) = k$, then, when depicting $f$ with a diagram, replace the node (or fill in the cell) corresponding to $p$ with $k$.

\begin{example}\label{labelling 1} If $(P, \leq)$ is the poset depicted by the left-most diagram shown below, where $P = \{p_1,p_2,p_3,p_4\}$, and if $f : P \to \mathbb{Z}$ is a function such that $f(p_1) = 4$, $f(p_2) = 2$, $f(p_3) = -1$, and $f(p_4) = 3$, then we will depict $f$ with the right-most diagram shown below.
\begin{center}
	\begin{tikzpicture}
	\node (d) at (0,0) {$p_4$};
	\node (a) at (-1,-1) {$p_1$};
	\node (b) at (0,-1) {$p_2$};
	\node (c) at (1,-1) {$p_3$};
	
	\draw  (a) -- (d) -- (c) -- (d) -- (b);
	\end{tikzpicture} $\;\;\;\;\;\;\;\;\;\;\;\;\;\;\;\;$
	\begin{tikzpicture}
	\node (d) at (0,0) {$3$};
	\node (a) at (-1,-1) {$4$};
	\node (b) at (0,-1) {$2$};
	\node (c) at (1,-1) {$-1$};
	
	\draw  (a) -- (d) -- (c) -- (d) -- (b);
	\end{tikzpicture}
\end{center}
\end{example}
\begin{example}\label{labelling 2} If $(P, \leq)$ is the poset depicted by the left-most diagram shown below, where $P = \{p_1,p_2,p_3\}$, and if $f : P \to \{1,2,3\}$ is a function such that $f(p_1) = 3$, $f(p_2) = 2$, and $f(p_3) = 2$, then we will depict $f$ with the right-most diagram shown below.
\begin{center}
	$\tableau{& p_1 \\ p_2 & p_3} \;\;\;\;\;\;\;\;\;\;\;\;\;\;\;\; \tableau{& 3 \\ 2 & 2}$
\end{center}
\end{example}

If $P$ is a finite poset and if $\omega$ is a labelling of $P$, then we would like to say that two elements $U_1, U_2 \in \mathcal{A}(P, \omega)$ are the same if the relative orderings of their entries are the same. Hence, we define the following equivalence relation.

\begin{definition}\label{order equivalence} Let $P$ be a poset, let $S_1$ and $S_2$ be subsets of $\mathbb{Z}$, and let $f_1 : P \to S_1$ and $f_2 : P \to S_2$ be maps. Then \emph{$f_1$ is order equivalent to $f_2$} if there exists an order-isomorphism $g : f_1(P) \to f_2(P)$ such that $f_2 = g \circ f_1$. Lastly, we write $f_1 \equiv f_2$ if $f_1$ is order equivalent to $f_2$, and write $f_1 \not\equiv f_2$ otherwise.
\end{definition}

\begin{example} Let $f_1 : P \to \mathbb{Z}$ and $f_2 : P \to  \mathbb{Z}$ be depicted by
$$\tableau{5 & 2 & 2 \\ & 6 & 7} \;\;\;\;\;\;\;\;\;\;\;\; \text{ and } \;\;\;\;\;\;\;\;\;\;\;\; \tableau{8 & -1 & -1 \\ & 9 & 15} $$
respectively. Then $f_1(P) = \{2,5,6,7\}$, $f_2(P) = \{-1,8,9,15\}$, and $g : f_1(P) \to f_2(P)$ is defined by $g(2) = -1$, $g(5) = 8$, $g(6) = 9$, and $g(7) = 15$. Since $f_2 = g \circ f_1$, and since $g$ is an order-isomorphism, $f_1$ is order equivalent to $f_2$.
\end{example}

We now formally define the enumeration problem for $(P, \omega)$-partitions that we are interested in analysing. Specifically, we are interested in enumerating the number $|Q, \omega|$ of equivalence classes as described below. For any set $S$, let $|S|$ denote the cardinality of $S$.

\begin{definition}\label{tableaux} Let $P$ and $Q$ be posets such that $Q$ is finite and $Q \subseteq P$. Moreover, let $\omega$ be a labelling $P$ and let $\equiv_{Q,\omega}$ denote the equivalence relation on $\mathcal{A}(Q, \omega)$ defined by $f \equiv_{Q,\omega} g$ if $f$ is order equivalent to $g$. Then define
$$\Tb(Q,\omega) = \mathcal{A}(Q, \omega)/\equiv_{Q, \omega}$$
and define
$$|Q, \omega| = |\Tb(Q, \omega)|. $$

Moreover, if $Q_1 \subseteq Q$ and if $T \in \Tb(Q, \omega)$, then let $T|_{Q_1}$ be the element $T'$ of $\Tb(Q_1, \omega)$ such that, for all $U \in \mathcal{A}(Q, \omega)$ satisfying $U \in T$, $U|_{Q_1} \in T'$.
\end{definition}

To depict an element $T$ of $\Tb(Q, \omega)$, where $Q$ and $\omega$ are as specified in Definition \ref{tableaux}, pick the element $U \in \mathcal{A}(Q, \omega)$ such that $U \in T$ and the range of $U : Q \to \mathbb{N}_0$ is $[m]$ for some integer $m \geq 1$, display the diagram that represents the element $U$, then say that the resulting diagram depicts $T$.

\begin{example}\label{tableaux example} Let $P$ be the six-element poset depicted by the left-most diagram shown below, let $Q$ be the five-element subposet of $P$ that is depicted by the cells in the left-most diagram that are filled with bullets, and let $\omega : P \to \{1,2,3,4,5,6\}$ be depicted by the right-most diagram shown below.
$$\tableau{\bullet & \bullet & \bullet \\ & \bullet & \bullet & {}} \;\;\;\;\;\;\;\;\;\;\;\;\;\;\;\; \tableau{2 & 1 & 3 \\ & 6 & 4 & 5} $$
For all $1 \leq k \leq 6$, let $p_k \in P$ satisfy $\omega(p_k) = k$. Then $\mathcal{A}(Q, \omega)$ consists of the order-preserving maps $f : Q \to \mathbb{N}_0$ such that $f(p_2) > f(p_1)$, $f(p_1) \geq f(p_3)$, $f(p_1) \geq f(p_6)$, $f(p_3) \geq f(p_4)$, and $f(p_6) > f(p_4)$. Three of the elements in $\mathcal{A}(Q, \omega)$ are depicted as follows.
$$U_1 \; = \; \tableau{3 & 2 & 2 \\ & 2 & 1 } \;\;\;\;\;\;\;\;\;\;\;\;\;\;\;\;\; U_2 \; = \; \tableau{5 & 4 & 3 \\ & 2 & 1} \;\;\;\;\;\;\;\;\;\;\;\;\;\;\;\;\; U_3 \; = \; \tableau{10 & 9 & 9 \\ & 9 & 7}$$
We have that $U_1 \equiv U_3$. However, $U_2 \not\equiv U_1$ and $U_2 \not\equiv U_3$. \\

Moreover, it can be checked that the $|Q, \omega| = 8$ elements of $\Tb(Q, \omega)$ are depicted by the below diagrams.
$$\tableau{3 & 2 & 2 \\ & 2 & 1 } \;\;\;\;\;\;\;\;\;\; \tableau{3 & 2 & 1 \\ & 2 & 1} \;\;\;\;\;\;\;\;\;\; \tableau{4 & 3 & 2 \\ & 3 & 1} \;\;\;\;\;\;\;\;\;\; \tableau{4 & 3 & 3 \\ & 2 & 1} $$ \newline
$$\tableau{4 & 3 & 2 \\ & 2 & 1 } \;\;\;\;\;\;\;\;\;\; \tableau{4 & 3 & 1 \\ & 2 & 1} \;\;\;\;\;\;\;\;\;\; \tableau{5 & 4 & 2 \\ & 3 & 1} \;\;\;\;\;\;\;\;\;\; \tableau{5 & 4 & 3 \\ & 2 & 1} $$
Now consider $P$. Note that $\Tb(P, \omega)$ is well-defined because $P$ is finite. So if $T$ is the element of $\Tb(P, \omega)$ depicted below
$$\tableau{4 & 3 & 3 \\ & 3 & 2 & 1} \;\; , $$
then $T|_Q$ is the element of $\Tb(Q, \omega)$ depicted below
 $$\tableau{3 & 2 & 2 \\ & 2 & 1} \;\;.$$
\end{example}

\section{Periodic P-Partitions}\label{sec:periodic}

In this section, we introduce a class of $(P, \omega)$-partitions that we call periodic $(P, \omega)$-partitions. Afterwards, we indicate how the problem of enumerating periodic $(P, \omega)$-partitions generalizes the problem of enumerating standard Young tableaux, and semi-standard Young tableaux, of parallelogramic shapes/shifted strips when the number of cells in each row is fixed. In addition to this, we describe examples that illustrate the level of generality that periodic $(P, \omega)$-partitions encompass. \\

We first introduce the notion of a connected triple.

\begin{definition}\label{connected triple} Let $P$ be a poset. Then a \emph{connected triple $(A,B,C)$ of $P$} is an ordered triple $(A,B,C)$ that satisfies the following two properties.
\begin{enumerate}
	\item The set $\{A, B, C\}$ is a set partition of the set of elements of $P$.

	\item For all $p \in A$ and for all $q \in C$, there exists an element $p' \in B$ such that $p < p' < q$ in $P$.
\end{enumerate}
Moreover, if $B$ is a non-empty subset of $P$, then $B$ \emph{connects} $P$ if there exist non-empty subsets $A$ and $C$ of $P$ such that $(A, B, C)$ is a connected triple of $P$.
\end{definition}

\begin{example}\label{connected triple example} A \emph{parallelogramic shape} \cite{CHC} is a Young diagram $X$ with $n$ rows and $k$ cells in each row such that for any two adjacent rows, the lower row is shifted one cell to the right of the upper row. For instance, if $n = 3$ and $k = 5$, then the corresponding parallelogramic shape $X$ is the following.
$$\tableau{{} & {} & {} & {} & {} \\ & {} & {} & {} & {} & {} \\ && {} & {} & {} & {} & {} } $$
Such shapes are also called \emph{shifted strips} where the parameter $k$ is called the \emph{width} of such a strip \cite{EYTPS}, and they were investigated in \cite{CHC, EYTPS, HAlgebra}. \\

We make the following observation which can be generalized to any parallelogramic shape. If $P$ is the poset corresponding to the parallelogramic shape with $k = 4$ and $n = 4$ depicted below, then if $A$ is the subset of $P$ that is depicted by the blank cells, if $B$ is the the subset of $P$ that is depicted by the cells filled with bullets, and if $C$ is the subset of $P$ that is depicted by the cells filled with asterisks, then $(A, B, C)$ is a connected triple of $P$.
$$\tableau{{} & {} & {} & {} \\ & {} & {} & {} & \bullet \\ && {} & \bullet & \bullet & \bullet \\ &&& * & * & * & * } $$

\end{example}

We extend the usual definition of the \emph{successor function} from the natural numbers to the integers. Namely, let $\s : \mathbb{Z} \to \mathbb{Z}$ be the successor function defined by $\s(n) = n+1$ for all integers $n$. With this function, we define periodic quadruple systems.

\begin{definition}\label{periodic quadruple system} Let $Z$ be a countably infinite poset, and let $\omega$ be a labelling of $Z$. Moreover, let $\pi : Z \to \mathbb{Z}$ be a surjective order preserving map, and let $\theta : Z \to Z$ be an order-automorphism on $Z$. Then the ordered quadruple $(Z, \omega, \pi, \theta)$ is a \emph{periodic quadruple system} if the following three properties hold.
\begin{enumerate}
\item There exists an order-automorphism $\alpha : \mathbb{Z} \to \mathbb{Z}$ such that the following diagram commutes.
\begin{center}
	\begin{tikzpicture}[scale = 2]
	\node (A) at (2,0) {$\mathbb{Z}$};
	\node (B) at (1,0) {$Z$};
	\node (C) at (0,0) {$\mathbb{Z}$};
	\node (D) at (2,-1) {$\mathbb{Z}$};
	\node (E) at (1,-1) {$Z$};
	\node (F) at (0,-1) {$\mathbb{Z}$};
	
	\draw [thick, ->] (B) -- (A) node [midway, above] {$\omega$};
	\draw [thick, ->] (B) -- (C) node [midway, above]{$\pi$};
	\draw [thick, ->] (E) -- (D) node [midway, below] {$\omega$};
	\draw [thick, ->] (E) -- (F) node [midway, below]{$\pi$};
	\draw [thick, ->] (A) -- (D) node [midway, right] {$\alpha$};
	\draw [thick, ->] (B) -- (E) node [midway, left] {$\theta$};
	\draw [thick, ->] (C) -- (F) node [midway, left] {$\s$};
	\end{tikzpicture}
\end{center}
\item For all integers $n$, $\pi^{-1}(\{n\})$ is finite.
\item There exists a finite subset $S$ of $Z$ such that $S$ connects $Z$. 
	
\end{enumerate}
	
\end{definition}

Definition \ref{periodic quadruple system} can be applied to (collections of) parallelogramic shapes.

\begin{example}\label{small 4 chain} Consider the poset $Z$ depicted by the left-most diagram in Figure \ref{fig:parallelogramic}. Moreover, let $\omega$ be the labelling of $Z$ depicted by the middle diagram in Figure \ref{fig:parallelogramic}. For all integers $n$, let $p_n$ denote the element of $Z$ such that $\omega(p_n) = n$. Next, define $\pi : Z \to \mathbb{Z}$ so that for all integers $n$, $\pi^{-1}(\{n\}) = \{p_{1-4n}, \, p_{2-4n}, \, p_{3-4n}, \, p_{4-4n} \}$. Furthermore, let $\theta : Z \to Z$ be the order-automorphism on $Z$ such that for all integers $n$, $\theta(p_n) = p_{n-4}$. \\

\begin{figure}
\begin{center}
	\begin{tikzpicture}[scale = 0.75]
	\node (a) at (0,0) {$\circ$};
	\node (b) at (1,1) {$\circ$};
	\node (c) at (2,2) {$\circ$};
	\node (d) at (3,3) {$\circ$};
	\node (e) at (0,-2) {$\circ$};
	\node (f) at (1,-1) {$\circ$};
	\node (g) at (2,0) {$\circ$};
	\node (h) at (3,1) {$\circ$};
	\node (i) at (0,-4) {$\circ$};
	\node (j) at (1,-3) {$\circ$};
	\node (k) at (2,-2) {$\circ$};
	\node (l) at (3,-1) {$\circ$};
	\node (m) at (0,-6) {$\circ$};
	\node (n) at (1,-5) {$\circ$};
	\node (o) at (2,-4) {$\circ$};
	\node (p) at (3,-3) {$\circ$};
	
	\node (bb) at (0.2,1.7) {$ $};
	\node (cc) at (1.2,2.7) {$ $};
	\node (dd) at (2.2,3.7) {$ $};
	\node (mm) at (0.7,-6.7) {$ $};
	\node (nn) at (1.7,-5.7) {$ $};
	\node (oo) at (2.7,-4.7) {$ $};
	
	\node (q) at (1.5,3) {$.$};
	\node (r) at (1.5,3.3) {$.$};
	\node (s) at (1.5,3.6) {$.$};
	
	\node (t) at (1.5,-6) {$.$};
	\node (u) at (1.5,-6.3) {$.$};
	\node (v) at (1.5,-6.6) {$.$};
	
	\draw (dd) -- (d) -- (c) -- (b) -- (a) -- (f) -- (k) -- (p) -- (o) -- (n) -- (m) -- (mm);
	
	\draw (cc) -- (c) -- (h) -- (g) -- (f) -- (e) -- (j) -- (o) -- (oo);
	
	\draw (bb) -- (b) -- (g) -- (l) -- (k) -- (j) -- (i) -- (n) -- (nn);
	
	\end{tikzpicture}
	\begin{tikzpicture}[scale = 0.75]
	\node (space) at (-3,0) {$ $};
	\node (a) at (0,0) {$0$};
	\node (b) at (1,1) {$-1$};
	\node (c) at (2,2) {$-2$};
	\node (d) at (3,3) {$-3$};
	\node (e) at (0,-2) {$4$};
	\node (f) at (1,-1) {$3$};
	\node (g) at (2,0) {$2$};
	\node (h) at (3,1) {$1$};
	\node (i) at (0,-4) {$8$};
	\node (j) at (1,-3) {$7$};
	\node (k) at (2,-2) {$6$};
	\node (l) at (3,-1) {$5$};
	\node (m) at (0,-6) {$12$};
	\node (n) at (1,-5) {$11$};
	\node (o) at (2,-4) {$10$};
	\node (p) at (3,-3) {$9$};
	
	\node (bb) at (0.2,1.7) {$ $};
	\node (cc) at (1.2,2.7) {$ $};
	\node (dd) at (2.2,3.7) {$ $};
	\node (mm) at (0.7,-6.7) {$ $};
	\node (nn) at (1.7,-5.7) {$ $};
	\node (oo) at (2.7,-4.7) {$ $};
	
	\node (q) at (1.5,3) {$.$};
	\node (r) at (1.5,3.3) {$.$};
	\node (s) at (1.5,3.6) {$.$};
	
	\node (t) at (1.5,-6) {$.$};
	\node (u) at (1.5,-6.3) {$.$};
	\node (v) at (1.5,-6.6) {$.$};
	
	\draw (dd) -- (d) -- (c) -- (b) -- (a) -- (f) -- (k) -- (p) -- (o) -- (n) -- (m) -- (mm);
	
	\draw (cc) -- (c) -- (h) -- (g) -- (f) -- (e) -- (j) -- (o) -- (oo);
	
	\draw (bb) -- (b) -- (g) -- (l) -- (k) -- (j) -- (i) -- (n) -- (nn);
	
	\end{tikzpicture}
	\begin{tikzpicture}[scale = 0.75]
	\node (space) at (-3,0) {$ $};
	\node (a) at (0,0) {$-6$};
	\node (b) at (1,1) {$-5$};
	\node (c) at (2,2) {$-4$};
	\node (d) at (3,3) {$-3$};
	\node (e) at (0,-2) {$-2$};
	\node (f) at (1,-1) {$-1$};
	\node (g) at (2,0) {$0$};
	\node (h) at (3,1) {$1$};
	\node (i) at (0,-4) {$2$};
	\node (j) at (1,-3) {$3$};
	\node (k) at (2,-2) {$4$};
	\node (l) at (3,-1) {$5$};
	\node (m) at (0,-6) {$6$};
	\node (n) at (1,-5) {$7$};
	\node (o) at (2,-4) {$8$};
	\node (p) at (3,-3) {$9$};
	
	\node (bb) at (0.2,1.7) {$ $};
	\node (cc) at (1.2,2.7) {$ $};
	\node (dd) at (2.2,3.7) {$ $};
	\node (mm) at (0.7,-6.7) {$ $};
	\node (nn) at (1.7,-5.7) {$ $};
	\node (oo) at (2.7,-4.7) {$ $};
	
	\node (q) at (1.5,3) {$.$};
	\node (r) at (1.5,3.3) {$.$};
	\node (s) at (1.5,3.6) {$.$};
	
	\node (t) at (1.5,-6) {$.$};
	\node (u) at (1.5,-6.3) {$.$};
	\node (v) at (1.5,-6.6) {$.$};
	
	\draw (dd) -- (d) -- (c) -- (b) -- (a) -- (f) -- (k) -- (p) -- (o) -- (n) -- (m) -- (mm);
	
	\draw (cc) -- (c) -- (h) -- (g) -- (f) -- (e) -- (j) -- (o) -- (oo);
	
	\draw (bb) -- (b) -- (g) -- (l) -- (k) -- (j) -- (i) -- (n) -- (nn);
	
	\end{tikzpicture}
\end{center}
\caption{}\label{fig:parallelogramic}
\end{figure}

To see that $(Z, \omega, \pi, \theta)$ satisfies Property 1 of Definition \ref{periodic quadruple system}, let $\alpha : \mathbb{Z} \to \mathbb{Z}$ be defined by $\alpha(n) = n - 4$ for all integers $n$. Then for all $p \in Z$, $\s(\pi(p)) = \pi(\theta(p))$ and $\alpha(\omega(p)) = \omega(\theta(p))$. For instance, $\s(\pi(p_7)) = \s(-1) = 0$, $\pi(\theta(p_7)) = 0$, $\alpha(\omega(p_7)) = \alpha(7) = 3$, and $\omega(\theta(p_7)) = \omega(p_3) = 3$, implying that $\s(\pi(p_7)) = 0 = \pi(\theta(p_7))$ and $\alpha(\omega(p_7)) = 3 = \omega(\theta(p_7))$. Moreover, by how $\pi$ is defined in this example, $|\pi^{-1}(\{n\})| = 4$ for all integers $n$. So $(Z, \omega, \pi, \theta)$ satisfies Property 2 of Definition \ref{periodic quadruple system}. To see that $(Z, \omega, \pi, \theta)$ satisfies Property 3 of Definition \ref{periodic quadruple system}, let $S \subset Z$ be defined by $S = \{p_1, p_2, p_3, p_5 \}$. The set $S$ is finite, and $S$ connects $Z$. Hence, $(Z, \omega, \pi, \theta)$ is a periodic quadruple system. \\

Lastly, let $\omega''$ be the labelling of $Z$ depicted by the right-most diagram in Figure \ref{fig:parallelogramic}. With $Z$, $\pi$, and $\theta$ being as described in this example, one can similarly check that $(Z, \omega'', \pi, \theta)$ is also a periodic quadruple system.
\end{example}

\begin{remark} The labelling $\omega''$ of $Z$ in Example \ref{small 4 chain} can be regarded as being a \emph{Schur labelling}. Roughly speaking, Schur labellings are labellings for \emph{skew shapes} that relate $(P, \omega)$-partitions to \emph{semi-standard Young tableaux of skew shapes} \cite{EC2}.
\end{remark}

Periodic quadruple systems can be very different from Example \ref{small 4 chain}. The following is an example of such a system.

\begin{figure}
	\begin{center}
		\begin{tikzpicture}[scale = 0.95]
		\node (A) at (0,0) {$\circ$};
		\node (B) at (-2,-2) {$\circ$};
		\node (C) at (1,-1) {$\circ$};
		\node (D) at (-1,-3) {$\circ$};
		\node (E) at (2,-2) {$\circ$};
		\node (F) at (0,-4) {$\circ$};
		\node (G) at (3,-3) {$\circ$};
		\node (GGG) at (3,-3.3) {$q_1$};
		\node (H) at (1,-5) {$\circ$};
		\node (HHH) at (1,-4.7) {$p_1$}; 
		\node (I) at (4,-4) {$\circ$};
		\node (III) at (4,-4.3) {$q_0$};
		\node (J) at (2,-6) {$\circ$};
		\node (JJJ) at (2,-5.7) {$p_0$};
		\node (K) at (5,-5) {$\circ$};
		\node (KKK) at (5,-5.3) {$q_{-1}$};
		\node (L) at (3,-7) {$\circ$};
		\node (LLL) at (3,-6.7) {$p_{-1}$};
		\node (M) at (6,-6) {$\circ$};
		\node (N) at (4,-8) {$\circ$};
		\node (O) at (7,-7) {$\circ$};
		\node (P) at (5,-9) {$\circ$};
		\node (Q) at (8,-8) {$\circ$};
		\node (R) at (6,-10) {$\circ$};
		
		\draw (B) arc (170:280:2.58);
		\draw (D) arc (170:280:2.58);
		\draw (F) arc (170:280:2.58);
		\draw (H) arc (170:280:2.58);
		\draw (J) arc (170:280:2.58);
		\draw (L) arc (170:280:2.58);
		
		\draw (N) arc (170:250:2.58);
		\draw (P) arc (170:218:2.58);
		\draw (R) arc (170:185:2.58);
		
		\draw (F) arc (280:200:2.58);
		\draw (D) arc (280:232:2.58);
		\draw (B) arc (280:265:2.58);
		
		\draw (A) arc (100:-10:2.58);
		\draw (C) arc (100:-10:2.58);
		\draw (E) arc (100:-10:2.58);
		\draw (G) arc (100:-10:2.58);
		\draw (I) arc (100:-10:2.58);
		\draw (K) arc (100:-10:2.58);
		
		\draw (M) arc (100:20:2.58);
		\draw (O) arc (100:52:2.58);
		\draw (Q) arc (100:85:2.58);
		
		\draw (E) arc (-10:70:2.58);
		\draw (C) arc (-10:38:2.58);
		\draw (A) arc (-10:5:2.58);
		
		\draw (B) -- (I);
		\draw (D) -- (K);
		\draw (F) -- (M);
		\draw (H) -- (O);
		\draw (J) -- (Q);
		
		\draw (A) -- (B);
		\draw (C) -- (D);
		\draw (E) -- (F);
		\draw (G) -- (H);
		\draw (I) -- (J);
		\draw (K) -- (L);
		\draw (M) -- (N);
		\draw (O) -- (P);
		\draw (Q) -- (R);
		
		\node (GG) at (-2.25,-1.25) {$ $};
		\draw (G) -- (GG);
		
		\node (EE) at (-1.75,-0.75) {$ $};
		\draw (E) -- (EE);
		
		\node (CC) at (-1.25,-0.25) {$ $};
		\draw (C) -- (CC);
		
		\node (AA) at (-0.75,0.25) {$ $};
		\draw (A) -- (AA);
		
		\node (LL) at (8.25,-8.75) {$ $};
		\draw (L) -- (LL);
		
		\node (NN) at (7.75,-9.25) {$ $};
		\draw (N) -- (NN);
		
		\node (PP) at (7.25,-9.75) {$ $};
		\draw (P) -- (PP);
		
		\node (RR) at (6.75,-10.25) {$ $};
		\draw (R) -- (RR);
		
		\node (aaa) at (-1.75,-0.25) {$.$};
		\node (bbb) at (-1.95,-0.05) {$.$};
		\node (ccc) at (-2.15,0.15) {$.$};
		
		\node (ddd) at (7.75,-9.75) {$.$};
		\node (eee) at (7.95,-9.95) {$.$};
		\node (fff) at (8.15,-10.15) {$.$};
		\end{tikzpicture}
	\end{center}
	\caption{A more exotic example of a periodic quadruple system.}\label{fig:exotic}
\end{figure}

\begin{example}\label{small exotic example} Let $Z$ be the poset depicted in Figure \ref{fig:exotic}. Moreover, let $p_1$, $q_1$, $p_0$, $q_0$, $p_{-1}$, and $q_{-1}$ be the six elements of $Z$ that are as specified in Figure \ref{fig:exotic}. Let $\omega : Z \to \mathbb{Z}$ be defined by $\omega(p_0) = 0$, $\omega(q_0) = 1$, and, for all $p \in Z$, $\omega(\theta(p)) = \omega(p) + 2$, let $\pi : Z \to \mathbb{Z}$ satisfy $\pi(p_1) = \pi(q_1) = 1$, $\pi(p_0) = \pi(q_0) = 0$, $\pi(p_{-1}) = \pi(q_{-1}) = -1$, and so on, and let $\theta : Z \to Z$ satisfy $\theta(p_{-1}) = p_0$, $\theta(q_{-1}) = q_0$, $\theta(p_0) = p_1$, $\theta(q_0) = q_1$, and so on. To check that the quadruple $(Z, \omega, \pi, \theta)$ satisfies Property 1 of Definition \ref{periodic quadruple system}, let $\alpha : \mathbb{Z} \to \mathbb{Z}$ be defined by $\alpha(n) = n + 2$ for all integers $n$. Then for all $p \in Z$, $\s(\pi(p)) = \pi(\theta(p))$ and $\alpha(\omega(p)) = \omega(\theta(p))$. \\

To see that $(Z, \omega, \pi, \theta)$ satisfies Property 3 of Definition \ref{periodic quadruple system}, let $P_n = \pi^{-1}([n])$ for all integers $n \geq 1$. We show that $P_{15}$ connects $Z$. Consider the element $q_0 \in Z$. We have the inequalities,
$$ q_0 < \theta^{4}(p_0) < \theta^{7}(p_0) < \theta^{10}(p_0) < \theta^{13}(p_0), $$
$$ q_0 < \theta^{4}(p_0) < \theta^{4}(q_0) < \theta^{8}(p_0) < \theta^{11}(p_0) < \theta^{14}(p_0),$$
and
$$ q_0 < \theta^{4}(p_0) < \theta^{4}(q_0) < \theta^{8}(p_0) < \theta^{8}(q_0) < \theta^{12}(p_0).$$

Moreover, $\{\theta^{14}(p_0), \theta^{13}(p_0), \theta^{12}(p_0) \}$ is the set of minimal elements of $\pi^{-1}([12,$ $ \infty))$. So as $\{\theta^{14}(p_0), \theta^{13}(p_0), \theta^{12}(p_0) \} \subset P_{15}$, it follows that for all $p \in Z$ satisfying $\pi(p) \geq 16$, there is an element $q \in P_{15}$ such that $q_0 < q < p$. Since $\theta$ is an order-automorphism on $Z$, the same conclusions also hold for $\theta^{-1}(q_0)$ and $\theta^{-1}(P_{15}) = \pi^{-1}([0, 14])$, and for $\theta^{-2}(q_0)$ and $\theta^{-2}(P_{15}) = \pi^{-1}([-1,13])$. So as $\{q_0, \theta^{-1}(q_0), \theta^{-2}(q_0) \}$ is the set of maximal elements of $\pi^{-1}(\mathbb{Z} \backslash \mathbb{N})$, it follows that $P_{15}$ connects $Z$. Hence, $(Z, \omega, \pi, \theta)$ is a periodic quadruple system.
\end{example}

We now introduce the $P$-partitions that we are interested in enumerating in this paper.

\begin{definition}\label{periodic P partition} Let $(Z, \omega', \pi, \theta)$ be a periodic quadruple system and let $n \in \mathbb{N}$. Then a \emph{length $n$ periodic $(P,\omega)$-partition derived from $(Z, \omega', \pi, \theta)$} is a $(P,\omega)$-partition $(P, \omega)$ such that $P = \pi^{-1}([n])$ and $\omega \equiv \omega'|_P$. Moreover, a \emph{periodic $(P, \omega)$-partition derived from $(Z, \omega', \pi, \theta)$} is a length $n$ periodic $(P, \omega)$-partition derived from $(Z, \omega', \pi, \theta)$ for some integer $n \geq 1$. Lastly, a \emph{periodic $(P, \omega)$-partition} is a periodic $(P, \omega)$-partition derived from $(Z, \omega', \pi, \theta)$ for some periodic quadruple system $(Z, \omega', \pi, \theta)$. For brevity, we also call periodic $(P, \omega)$-partitions \emph{periodic $P$-partitions}.
\end{definition}

\begin{remark} Since $P = \pi^{-1}([n])$, $P$ can be informally thought of as the poset that results from pasting together $n$ copies $\pi^{-1}(\{1\})$, $\pi^{-1}(\{2\})$, $\dots$, and $\pi^{-1}(\{n\})$ of the poset $\pi^{-1}(\{{0}\})$ where the pasting is determined by the periodic quadruple system $(Z, \omega, \pi, \theta)$.
\end{remark}

Enumeration formulas for counting standard Young tableaux on certain non-classical shapes was established in \cite{CHC} with special cases being established in \cite{EYTPS, HAlgebra}. We describe these objects below. Let $X$ be a parallelogramic shape with $n$ rows and $k$ cells in each row as described in Example \ref{connected triple example}. Then a \emph{standard Young tableau of shape $X$} is a bijective filling of the cells of $X$ with integers from $[nk]$ such that the entries increase along every row of $X$ from left to right and the entries increase along every column of $X$ from top to bottom, and we call such a standard Young tableau a \emph{standard Young tableau of parallelogramic shape}. For instance, two standard Young tableaux of parallelogramic shape are displayed below.
$$\tableau{1 & 2 & 4 & 6 \\ & 3 & 5 & 7 & 10 \\ && 8 & 9 & 11 & 12} \;\;\;\;\;\;\;\;\;\;  \tableau{1 & 2 & 4 & 5 \\ & 3 & 6 & 8 & 9 \\ && 7 & 10 & 11 & 12}$$
The enumeration problem we are considering in this paper generalizes the problem of enumerating standard Young tableaux of parallelogramic shape when $k$ is fixed.

\begin{example}\label{standard periodic ex} Let $Z$, $\omega$, $\pi$, and $\theta$ be as in Example \ref{small 4 chain}. Then the length $n$ periodic $(P, \omega)$-partitions derived from $(Z, \omega, \pi, \theta)$ correspond to the standard Young tableaux of parallelogramic shape with $n$ rows and four cells in each row. Specifically, if $n \geq 1$ is an integer and if $P_n = \pi^{-1}([n])$, then the elements of $\Tb(P_n, \omega)$ correspond to the aforementioned standard Young tableaux. For instance, when $n = 3$, three elements of $\Tb(P_3, \omega)$ are depicted below.
\begin{center}
$\tableau{12 & 11 & 9 & 8 \\ & 10 & 7 & 6 & 5 \\ && 4 & 3 & 2 & 1} \;\;\;\;\; \tableau{12 & 11 & 9 & 8 \\ & 10 & 7 & 6 & 5 \\ && 4 & 3 & 2 & 1} \;\;\;\;\; \tableau{12 & 11 & 10 & 8 \\ & 9 & 7 & 5 & 4 \\ && 6 & 3 & 2 & 1}$
\end{center}
In particular, replacing every entry $m$ in each of the three diagrams depicted above with $12 - m + 1$ gives three examples of standard Young tableaux of paralellogramic shape with three rows and four cells in each row.
\end{example}

Next, we explain how we can also apply our results to semi-standard Young tableaux. If $X$ is a parallelogramic shape, then define a \emph{semi-standard Young tableau of shape $X$} to be a function $f : X \to \mathbb{Z}$ where $f(i_1, j) < f(i_2,j)$ if $(i_1,j), (i_2,j) \in X$ and $i_1 < i_2$, and where $f(i,j_1) \leq f(i,j_2)$ if $(i,j_1), (i,j_2) \in X$ and $j_1 < j_2$. That is, fill the cells of $X$ so that the entries weakly increase along every row of $X$ from left to right, and the entries strictly increase along every column of $X$ from top to bottom. Moreover, if $f$ is a semi-standard Young tableau of shape $X$, then define a \emph{semi-standard tableau class on $X$} to be the set $F$ of semi-standard Young tableau of shape $X$ that are order equivalent to $f$. For example, if $X$ is a parallelogramic shape with $n = 3$ and $k = 4$, then three semi-standard tableau classes on $X$ are depicted below.
$$\tableau{1 & 1 & 2 & 2 \\ & 2 & 3 & 4 & 4 \\ && 4 & 5 & 5 & 5} \;\;\;\;\;\;\; \tableau{1 & 1 & 1 & 1 \\ & 2 & 2 & 2 & 2 \\ && 3 & 3 & 3 & 3} \;\;\;\;\;\; \tableau{1 & 2 & 4 & 6 \\ & 3 & 5 & 7 & 10 \\ && 8 & 9 & 11 & 12}$$
The periodic $P$-partitions that we introduced in Definition \ref{periodic P partition} can be regarded as a generalization of the semi-standard Young tableau of parallelogramic shape $X$ if the number of cells in each row of $X$ is fixed. This is explained in the following example when the number of cells in each row is four.

\begin{example}\label{semi-standard tableau periodic quadruple} Let $Z$, $\omega''$, $\pi$, and $\theta$ be as in Example \ref{small 4 chain}. Then the length $n$ periodic $(P, \omega)$-partitions derived from $(Z, \omega, \pi, \theta)$ correspond to the semi-standard tableau classes on $X$ when $X$ is a parallelogramic shape with $n$ rows and four cells in each row. For instance, when $n = 3$, three elements of $\Tb(P_3, \omega'')$ are depicted below.
\begin{center}
$\tableau{6 & 6 & 6 & 6 \\ & 5 & 4 & 3 & 2 \\ && 3 & 2 & 1 & 1} \;\;\;\;\; \tableau{5 & 5 & 5 & 5 \\ & 4 & 3 & 2 & 2 \\ && 2 & 1 & 1 & 1} \;\;\;\;\; \tableau{5 & 5 & 4 & 4 \\ & 3 & 2 & 2 & 2 \\ && 1 & 1 & 1 & 1}$
\end{center}
In particular, replacing every entry $k$ in each of the three diagrams depicted above with $12 - k + 1$ gives three examples of semi-standard tableaux classes of paralellogramic shape with three rows and four cells in each row.
\end{example}

In the last part of this section, we describe how more examples of periodic quadruple systems can be generated if they are subposets of $\mathbb{Z}^d$.

\begin{definition}\label{Euclidean PQS} Fix a positive integer $d \geq 2$ and consider $\mathbb{Z}^d$. Next, let $X$ be a finite subset of $\mathbb{Z}^d$. Moreover, let $\Delta = (k_1, k_2, \dots, k_d) \in \mathbb{N}^d$ be such that $(X + \Delta) \cap X = \emptyset$. Then define
$$Z(X, \Delta) = \bigcup_{n \in \mathbb{Z}} (X + n \, \Delta).$$
\end{definition}

Informally, $Z(X, \Delta)$ is a pairwise-disjoint union of translates of a finite subset $X$ of $\mathbb{Z}^d$ where $\Delta$ is the translate. An important property of this translate is that all coordinates of $\Delta$ are positive.

\begin{remark} Note that by our convention, $Z(X, \Delta)$ is a subposet of $\mathbb{Z}^d$, where the partial order on $\mathbb{Z}^d$ is as described in Section \ref{sec:prelims}, because $Z(X, \Delta)$ is a subset of $\mathbb{Z}^d$.
\end{remark}

\begin{example}\label{example for ZXDelta} If $d = 2$, $X = \{(1,0), (2,0), (3,0), (4,0)\}$, and $\Delta = (1,1)$, then the subposet $Z(X, \Delta)$ of $\mathbb{Z}^2$ is the poset depicted in the left-most diagram referred in Example \ref{small 4 chain}.
\end{example}

The posets $Z(X, \Delta)$ satisfy the following property.

\begin{lemma}\label{Euclidean System} Fix a positive integer $d \geq 2$ and let $Z(X, \Delta)$ be as described in Definition \ref{Euclidean PQS}. Then there exists a finite set $S$ that connects $Z(X, \Delta)$.
\end{lemma}

\begin{example}\label{ex:small 4 chain as Z} The poset $Z$ in Example \ref{small 4 chain} is an example of a poset $Z(X, \Delta)$ where $X = \{(1,0), (2,0), (3,0), (4,0)\}$ and $\Delta = (1,1)$.
\end{example}

\begin{proof} Define $\pi : Z(X, \Delta) \to \mathbb{Z}$ by $\pi(p + n \Delta) = n$ for all $p \in X$ and for all $n \in \mathbb{Z}$. By Definition \ref{Euclidean PQS}, $(X + \Delta) \cap X = \emptyset$. So $\pi$ is well-defined. Hence, we proceed as follows. Because all coordinates of $\Delta$ are positive, there is a sufficiently large integer $k \geq 1$ such that for all $p,q \in X$, the inequality
\begin{equation}\label{shift inequality}
 p \leq q + k \Delta 
\end{equation}
holds in $Z(X, \Delta)$. In particular, for all $p,q \in Z(X, \Delta)$ satisfying $\pi(q) \geq \pi(p) + k$, Inequality \ref{shift inequality} implies that $p < q$ in $Z(X, \Delta)$. Let
$$S = \bigcup_{n=1}^{2k+1} (X + n \Delta).$$
Since $X$ and $k$ are finite, $S$ is finite. For all $p,q \in Z(X,\Delta)$ satisfying $\pi(p) < 1$ and $\pi(q) > 2k+1$, the following holds. Let $p' \in Z(X,\Delta)$ satisfy $\pi(p') = k+1$. Since $\pi(p') \geq \pi(p) + k$ and $\pi(q) \geq \pi(p') + k$, $p < p' < q$ in $Z$. Hence, by Definition \ref{connected triple}, $S$ connects $Z(X, \Delta)$.
\end{proof}

There is at least one periodic quadruple system that can be constructed from the poset $Z(X, \Delta)$ in Lemma \ref{Euclidean System}.

\begin{proposition}\label{Full Euclidean system} Fix a positive integer $d \geq 2$, and let $Z(X, \Delta)$ be as described in Definition \ref{Euclidean PQS}. Moreover, let $\omega_X$ be a dual-natural labelling of $X$, let $\pi : Z(X, \Delta) \to \mathbb{Z}$ be defined by $\pi(p + n\Delta) = n$ for all $p \in X$ and for all $n \in \mathbb{Z}$, and let $\theta : Z(X, \Delta) \to Z(X, \Delta)$ be defined by $\theta(p) = p + \Delta$ for all $p \in Z(X, \Delta)$. Then there exists a dual-natural labelling $\omega : Z(X, \Delta) \to \mathbb{Z}$ of $Z(X, \Delta)$ such that $\omega_X = \omega|_X$ and $(Z(X, \Delta), \omega, \pi, \theta)$ is a periodic quadruple system.
\end{proposition}

\begin{example} The periodic quadruple system $(Z, \omega, \pi, \theta)$ in Example \ref{small 4 chain} is an example of a periodic quadruple system constructed in Proposition \ref{Full Euclidean system}.
\end{example}

\begin{proof} By Definition \ref{Euclidean PQS}, $(X + \Delta) \cap X = \emptyset$. So $\pi$ is well-defined. Moreover, note that $\pi$ is order-preserving. So we proceed as follows. For all $p \in Z(X, \Delta)$, $\pi(\theta(p)) = \pi(p + \Delta) = \pi(p) + 1$. So $\pi \circ \theta = \s \circ \pi$ where $\s$ is as described in Definition \ref{periodic quadruple system}. Note that $\theta$ is an order-automorphism on $Z(X, \Delta)$. Hence, we define the following. Let $\omega : Z(X, \Delta) \to \mathbb{Z}$ be the bijection such that $\omega|_X = \omega_X$, $\omega(p) > \omega(q)$ for all $p,q \in Z(X, \Delta)$ such that $\pi(p) < \pi(q)$, and $(\omega \circ \theta^n)|_X \equiv \omega|_X$ for all $n \in \mathbb{Z}$. As $(\omega \circ \theta^n)|_X \equiv \omega_P$ for all $n \in \mathbb{Z}$, $\omega|_{\theta^n(X)}$ is order-reversing for all $n \in \mathbb{Z}$. So if $\pi(p) = \pi(q)$, then $\omega(p) > \omega(q)$ because $p < q$ in $Z(X, \Delta)$. If $\pi(p) \neq \pi(q)$, then $\pi(p) < \pi(q)$ since $\pi$ is order-preserving and since $p < q$. It follows from the definition of $\omega$ that $\omega(p) > \omega(q)$. So $\omega$ is dual-natural. \\

Because $(\omega \circ \theta^n)|_X \equiv \omega|_X$ for all $n \in \mathbb{Z}$ and because $\omega$ is a bijection, it follows that for all $p \in Z(X, \Delta)$, $\omega(p + \Delta) = \omega(p) - |X|$. So for all $p \in Z(X, \Delta)$, $\omega(\theta(p)) = \omega(p + \Delta) = \omega(p) - |X|$. Hence, the map $\alpha : \mathbb{Z} \to \mathbb{Z}$ defined by $\alpha(n) = n - |X|$ satisfies $\omega \circ \theta = \alpha \circ \omega$. It follows that the quadruple $(Z, \omega, \pi, \theta)$ satisfies Property 1 of Definition \ref{periodic quadruple system}. For all $n \in \mathbb{Z}$, $|\pi^{-1}(\{n\})| = |X + n \Delta| = |X|$, implying that $|\pi^{-1}(n)|$ is finite. So $(Z, \omega, \pi, \theta)$ satisfies Property 2 of Definition \ref{periodic quadruple system}. Lastly, Lemma \ref{Euclidean System} implies that $(Z, \omega, \pi, \theta)$ satisfies Property 3 of Definition \ref{periodic quadruple system}. From this, the lemma follows.
\end{proof}

Lastly, we make the following remarks.

\begin{remark} Given any periodic quadruple system $(Z, \omega, \pi, \theta)$, the quadruple $(Z, \omega', \pi, \theta)$, where $\omega'(p) = -\omega(p)$ for all $p \in Z$, is also a periodic quadruple system. In particular, dual versions of the periodic quadruple systems in Example \ref{small 4 chain}, Example \ref{small exotic example}, and Proposition \ref{Full Euclidean system} also exist.
\end{remark}

\begin{remark} Standard Young tableaux on \emph{truncated shifted shapes} have recently been a topic of interest \cite{AKR, SYT, ARTHS, PTT}. Using Proposition \ref{Full Euclidean system}, it can be checked that by defining $X$ and $\Delta$ appropriately, the results in this paper can be applied to Standard Young tableaux on truncated shifted shapes in which the numbers of cells in each row, when read from top to bottom, form (part of) a periodic sequence.
\end{remark}

\section{Matrix Difference Equations and Recurrence Relations}\label{sec:matrix}

In this section, we construct a matrix with non-negative integer entries that we call a \emph{tableau transfer matrix}, then prove that the periodic $(P, \omega)$-partitions can be enumerated with these matrices. We then explore some of the implications of our result by giving a new proof of L\'opez et.al. and Sun's constant coefficient linear recurrence results and by explaining how to derive semi-standard variants of their results.

\subsection{The Construction}

Building from the previous section, we define the following.

\begin{definition}\label{index shapes} Let $(Z, \omega, \pi, \theta)$ be a periodic quadruple system. Then an \emph{index shape of $(Z, \omega, \pi, \theta)$} is a finite subset $S$ of $\pi^{-1}(\mathbb{Z} \backslash \mathbb{N})$ that satisfies the following two properties.
\begin{enumerate}
	\item $(\pi^{-1}(\mathbb{Z} \backslash \mathbb{N}) \, \backslash \, S \,,\; S\,,\; \pi^{-1}(\mathbb{N}) )$ is a connected triple of $Z$.
	\item $\theta(S) \subseteq S \cup \pi^{-1}(\{1\})$.
\end{enumerate}
When $(Z, \omega, \pi, \theta)$ is not specified, we will call $S$ an \emph{index shape}.
\end{definition}

\begin{remark}\label{rem:index shapes always exist} Index shapes always exist. Let $(Z, \omega, \pi, \theta)$ be a periodic quadruple system. Then by Definition \ref{periodic quadruple system}, there exists a finite set $S$ that connects $Z$. So as $\theta$ is an order-automorphism on $Z$, it follows that for some integer $n$, $\theta^n(S)$ connects $Z$ and $\theta^n(S) \subseteq \pi^{-1}([n_1,n_2])$ for some integers $n_1 \leq n_2 \leq 0$. It can be checked that $\pi^{-1}([n_1,n_2])$ is an index shape of $(Z, \omega, \pi, \theta)$.
\end{remark}

\begin{remark} If $S$ is an index shape of a periodic quadruple system $(Z, \omega, \pi, \theta)$, then $S$ connects $Z$. This fact, and the fact that $S$ is finite, will become important later on in this section.
\end{remark}

\begin{example}\label{index shape ex} Let $(Z, \omega, \pi, \theta)$ and $S$ be as defined in Example \ref{small 4 chain}. Moreover, for all integers $n$, let $p_n \in Z$ satisfy $\omega(p_n) = n$. In Example \ref{small 4 chain}, we saw that $\theta(p_n) = p_{n-4}$, that $S = \{p_1, p_2, p_3, p_5 \}$, and that $(\pi^{-1}(\mathbb{Z} \backslash \mathbb{N}) \, \backslash \, S , S , \pi^{-1}(\mathbb{N}))$ is a connected triple of $Z$. Hence, $S$ satisfies Property 1 of Definition \ref{index shapes}. Since $\theta(S) = \{p_{-3}, p_{-2}, p_{-1}, p_1\}$ and $\pi^{-1}(\{1\}) = \{p_{-3}, p_{-2}, p_{-1}, p_0 \}$, 
$$\theta(S) = \{p_{-3}, p_{-2}, p_{-1}, p_1\} \subseteq \{p_{-3}, p_{-2}, p_{-1}, p_0, p_1, p_2, p_3, p_5 \} = S \cup \pi^{-1}(\{1\}).$$
Hence, $\theta(S) \subseteq S \cup \pi^{-1}(\{1\})$, implying that $S$ satisfies Property 2 of Definition \ref{index shapes}. It follows that $S$ is an index shape of $(Z, \omega, \pi, \theta)$.
\end{example}

The following definition is an analogue of the notion of order equivalence that depends on the order-automorphism $\theta$ on $Z$ in a periodic quadruple system.

\begin{definition}\label{similarity} Let $(Z, \omega, \pi, \theta)$ be a periodic quadruple system, and let $S_1$ and $S_2$ be subsets of $Z$. Moreover, assume that $S_2 = \theta^k(S_1)$ for some integer $k$. Then for all $T_1 \in \Tb(S_1, \omega)$ and $T_2 \in \Tb(S_2, \omega)$, write $T_1 \equiv_\theta T_2$ if for all $U_1 \in T_1$ and $U_2 \in T_2$, there is an order-isomorphism $g : U_1(S_1) \to U_2(S_2)$ such that $g \circ U_1 = U_2 \circ \theta^k$.
\end{definition}

Informally, the above definition says that $T_1 \equiv_\theta T_2$ if the relative ordering of the entries in $T_1$ is the same as the relative ordering of the entries in $T_2$.

\begin{example}\label{illustrative ex} Let $(Z, \omega, \pi, \theta)$ and $S$ be as described in Example \ref{index shape ex}. The eight-element subposet $S \, \cup \, \theta^{-1}(S) \, \cup \, \pi^{-1}(\{0\})$ of $Z$ is depicted by the Young diagram below. In the below diagram, $S$ is depicted by the cells that are filled with circles or asterisks, and $\theta^{-1}(S)$ is depicted by the cells that are filled with asterisks and bullets.
\begin{center}
$\tableau{& \circ \\ \circ & \circ & * \\ {} & \bullet & \bullet & \bullet }$
\end{center}
Next, let $T$ be the element of $\Tb(S \, \cup \, \theta^{-1}(S) \, \cup \, \pi^{-1}(\{0\}) , \omega)$ that is depicted by the diagram shown below. Lastly, for all integers $1 \leq i \leq 8$, let $c_i$ denote the element of $S \, \cup \, \theta^{-1}(S) \, \cup \, \pi^{-1}(\{0\})$ that is depicted by the cell in the diagram below that contains the entry $i$. 
\begin{center}
	$\tableau{& 7 \\ 8 & 6 & 3 \\ 5 & 4 & 2 & 1 }$
\end{center}
The function $\theta$ satisfies $\theta(c_7) = c_3$, $\theta(c_8) = c_4$, $\theta(c_6) = c_2$, and $\theta(c_3) = c_1$. From this, it follows that $T|_{\theta^{-1}(S)} \equiv_\theta T|_S$. For instance, let $U_1 \in T|_{\theta^{-1}(S)}$ be defined by $U_1(c_7) = 7$, $U_1(c_8) = 8$, $U_1(c_6) = 6$, and $U_1(c_3) = 3$, and let $U_2 \in T|_S$ be defined by $U_2(c_3) = 3$, $U_2(c_4) = 4$, $U_2(c_2) = 2$, and $U_2(c_1) = 1$. Then the order-isomorphism $g : U_1(\theta^{-1}(S)) \to U_2(S)$ defined by $g(3) = 1$, $g(6) = 2$, $g(7) = 3$, and $g(8) = 4$ satisfies $U_2 \circ \theta^k = g \circ U_1$ with $k = 1$.
\end{example}

We now use index shapes to define the following family of square matrices. Informally, we are defining matrices that are built from index shapes and that allow us to enumerate many different $P$-partitions at once. Recall that we write $M(i,j)$ to denote the entry in row $i$ and column $j$ of a matrix $M$. Lastly, in the following definition, note that $S \subseteq \theta^{-1}(S) \cup \pi^{-1}(\{0\})$ due to Property 2 of Definition \ref{index shapes} and the fact that $\theta$ is an order-automorphism on $Z$.

\begin{definition}\label{tableau transfer matrix} Let $(Z, \omega, \pi, \theta)$ be a periodic quadruple system, and let $S$ be an index shape of $(Z, \omega, \pi, \theta)$. Moreover, let $L$ be an indexing of $\Tb(S, \omega)$, let $N = |S, \omega|$, and let $R : [N] \to \Tb(S, \omega)$ be the inverse of $L$. Then the \emph{tableau transfer matrix $M$ derived from $(Z, \omega, \pi, \theta)$, $S$, and $L$} is the $N$ by $N$ matrix $M$ such that for all $1 \leq i \leq N$ and $1 \leq j \leq N$, $M(i,j)$ is the number of elements $T$ of $$\Tb(\theta^{-1}(S) \cup \pi^{-1}(\{0\}), \omega)$$ such that 
\begin{center}	
$T|_{S} \equiv_\theta R(i)$ \; and \; $T|_{\theta^{-1}(S)} \equiv_\theta R(j)$.
\end{center}

When $L$ is not specified, we will call $M$ a \emph{tableau transfer matrix derived from $(Z, \omega, \pi, \theta)$ and $S$}, when $S$ and $L$ are not specified, we will call $M$ a \emph{tableau transfer matrix derived from $(Z, \omega, \pi, \theta)$}, and when $(Z, \omega, \pi, \theta)$, $S$, and $L$ are not specified, we will call $M$ a \emph{tableau transfer matrix}.
\end{definition}

\begin{example}\label{matrix example} Consider the periodic quadruple system $(Z, \omega, \pi, \theta)$ and the index shape $S$ from Example \ref{illustrative ex}. There are two elements of $\Tb(S, \omega)$, and they are depicted by the below diagrams.
\begin{center}
	$\tableau{& 4 \\ 3 & 2 & 1}\;\;\;\;\;\;\;\;\;\;$
	$\tableau{& 3 \\ 4 & 2 & 1}$
\end{center}

In particular, $N = |S, \omega| = 2$. Next, define $L : \Tb(S, \omega) \to \{1, 2\}$ so that $L$ sends the element of $\Tb(S, \omega)$ depicted by the left-most diagram shown above to $1$ and $L$ sends the element of $\Tb(S, \omega)$ depicted by the right-most diagram shown above to $2$. Lastly, let $R = L^{-1}$ be the inverse of $L$ so that
\begin{center}
	$R(1) \;\; = \;\; \tableau{& 4 \\ 3 & 2 & 1} \;\;\;\;\;\;\; \text{ and } \;\;\;\;\;\; R(2) \;\; = \;\; \tableau{& 3 \\ 4 & 2 & 1} \;\;\; .$
\end{center}

The tableau transfer matrix $M$ derived from $(Z, \omega, \pi, \theta)$, $S$, and $L$ is equal to
$$\begin{bmatrix}
3 & 4 \\
2 & 3 \\
\end{bmatrix} \;. $$
To see how to construct $M$, we calculate $M(2,1)$. The integer $M(2,1)$ is the number of elements $T \in \Tb(S \, \cup \, \theta^{-1}(S) \, \cup \, \pi^{-1}(\{0\}) , \omega)$ such that $T|_S = R(2)$ and $T|_{\theta^{-1}(S)} \equiv_\theta R(1)$. By using the depictions of $S \, \cup \, \theta^{-1}(S) \, \cup \, \pi^{-1}(\{0\})$, $S$, and $\theta^{-1}(S)$ from Example \ref{illustrative ex}, it can be checked by an exhaustive search that the elements $T$ of $\Tb(S \, \cup \, \theta^{-1}(S) \, \cup \, \pi^{-1}(\{0\}), \omega)$ that satisfy $T|_S = R(2)$ and $T|_{\theta^{-1}(S)} \equiv_\theta R(1)$ are the two elements depicted by the below diagrams.

\begin{center}
	$\tableau{& 8 \\ 7 & 6 & 3 \\ 5 & 4 & 2 & 1 } \;\;\;\;\;\;$
	$\tableau{& 8 \\ 7 & 5 & 3 \\ 6 & 4 & 2 & 1 }$
\end{center}
	
\end{example}

We first prove a lemma that will be useful for the rest of this paper.

\begin{lemma}\label{bipartition lemma} Let $P$ be a finite poset and let $\omega_P$ be a labelling of $P$. Moreover, let $S_0$ be an order ideal of $P$ and let $S_1$ be an order filter of $P$ such that $\{S_0, S_1\}$ is a set partition of the set of elements of $P$. Lastly, let $\omega_0$ be a labelling of $S_0$ such that $\omega_0 \equiv \omega_P|_{S_0}$ and let $\omega_1$ be a labelling of $S_1$ such that $\omega_1 \equiv \omega_P|_{S_1}$. Then for all $T_0 \in \Tb(S_0, \omega_0)$ and for all $T_1 \in \Tb(S_1, \omega_1)$, there exists an element $T \in \Tb(S, \omega)$ such that $T|_{S_0} = T_0$, $T|_{S_1} = T_1$, and for all $U \in T$, $p \in S_0$, and $q \in S_1$, $U(p) > U(q)$.
\end{lemma}

\begin{example} Let $P$ be the eight-element poset depicted by the below Young diagram, let $S_0$ be depicted by the blank cells, and let $S_1$ be depicted by the cells that are filled with bullets.
\begin{center}
$\tableau{{} & {} & {} & {} \\ & \bullet & \bullet & \bullet & \bullet}$
\end{center}
Moreover, let $\omega$ be the following labelling of $P$
\begin{center}
$\tableau{1 & 2 & 3 & 4 \\ & 5 & 6 & 7 & 8}\;\;\;$.
\end{center}
Next, let $\omega_0$ be a labelling of $S_0$ such that $\omega_0 \equiv \omega|_{S_0}$ and let $\omega_1$ be a labelling of $S_1$ such that $\omega_1 \equiv \omega|_{S_1}$. Lastly, let $T_0 \in \Tb(S_0, \omega_0)$ be depicted by the left-most diagram below and let $T_1 \in \Tb(S_1, \omega_1)$ be depicted by the right-most diagram below.
\begin{center}
$\tableau{2 & 2 & 2 & 1}$ \;\;\;\;\;\;\; $\tableau{2 & 1 & 1 & 1}$
\end{center}
Then an element $T \in \Tb(S, \omega)$ as described in Lemma \ref{bipartition lemma} is depicted below.
\begin{center}
$\tableau{4 & 4 & 4 & 3 \\ & 2 & 1 & 1 & 1}$
\end{center}
\end{example}

\begin{proof} Let $T_0 \in \Tb(S_0, \omega_0)$ and let $T_1 \in \Tb(S_1, \omega_1)$.  Because $S_0$ and $S_1$ are finite and because $S_0 \cap S_1 = \emptyset$, there exist elements $U_0 \in T_0$ and $U_1 \in T_1$ such that for all $p \in S_0$ and $q \in S_1$, $U_0(p) > U_0(q) \geq 0$. So define the map $U : S \to \mathbb{N}_0$ by letting $U(p) = U_0(p)$ if $p \in S_0$ and $U(p) = U_1(p)$ if $p \in S_1$. This function $U$ satisfies $U|_{S_0} \equiv U_0$ and $U|_{S_1} \equiv U_1$. Because $S_0$ is an order ideal of $P$ and because $S_1$ is an order filter of $P$, the fact that $U(p) = U_0(p) > U_0(q) = U(q)$ for all $p \in S_0$ and $q \in S_1$ implies that $U \in \mathcal{A}(P, \omega_P)$. Let $T$ be the element of $\Tb(P, \omega_P)$ that contains $U$. Then $T$ satisfies the conditions described in the lemma.
\end{proof}

To complement Definition \ref{tableau transfer matrix}, we define the following. Informally, we are defining the column vectors that correspond to the tableau transfer matrices. Recall that we write $v(i)$ to denote the entry in row $i$ of a column vector $v$.

\begin{definition}\label{combinatorial vectors} Let $(Z, \omega, \pi, \theta)$ be a periodic quadruple system, let $S$ be an index shape of $(Z, \omega, \pi, \theta)$, and let $L$ be an indexing of $\Tb(S, \omega)$. Moreover, let $R = L^{-1}$ be the inverse of $L$. Then an \emph{admissible number for $(Z, \omega, \pi, \theta)$ and $S$} is an integer $n'$ such that $n' \leq 0$ and $S \subseteq \pi^{-1}([n',0])$. Next, let $n'$ be an admissible number for $(Z, \omega, \pi, \theta)$ and $S$, let $P' = \pi^{-1}([n',0])$, let $T_0 \in \Tb(P', \omega)$, let $n \in \mathbb{N}$, let $Q_n = \pi^{-1}([n',n])$, and let $P_n = \pi^{-1}([n])$. \\
	
Define \emph{the $n^{th}$ set derived from $(Z, \omega, \pi, \theta)$, $S$, and $T_0$} to be the set $X_n(T_0)$ of elements $T \in \Tb(Q_n, \omega)$ such that $T|_{P'} = T_0$ and the following condition holds. If $U \in T$, $p \in P'$, and $q \in P_n$, then $U(p) > U(q)$. \\

Moreover, for all $1 \leq i \leq |S, \omega|$, let the \emph{$i^{th}$ part of the $n^{th}$ set derived from $(Z, \omega, \pi, \theta)$, $S$, $L$, and $T_0$} be the set $X_{n,i}(T_0)$ of elements $T \in X_n(T_0)$ such that $T|_{\theta^n(S)} \equiv_\theta R(i)$. Lastly, define the \emph{$n^{th}$ vector derived from $(Z, \omega, \pi, \theta)$, $S$, $L$, and $T_0$} to be the column vector $v_n$ with $|S, \omega|$ entries such that, for all $1 \leq i \leq |S, \omega|$, $v_n(i) = |X_{n,i}(T_0)|$.

\end{definition}

Informally, the above definition describes the following. The $n^{th}$ sets as given in the above definition are a collection of modified periodic $P$-partitions that will allow us to enumerate the periodic $P$-partitions themselves. Moreover, the $i^{th}$ parts of such sets, as given in the above definition, provides us with a set partition of such collections of modified periodic $P$-partitions that we will use later on in this section, and the $n^{th}$ vectors give the cardinalities of these $i^{th}$ parts. Lastly, admissible numbers enable us to effectively use such modified periodic $P$-partitions.

\begin{example}\label{combinatorial vector example} Let $(Z, \omega, \pi, \theta)$, $S$, and $L$ be as in Example \ref{matrix example}, and let $R = L^{-1}$ be the inverse of $L$. Moreover, for all integers $n$, let $p_n \in Z$ satisfy $\omega(p_n) = n$. As described in Example \ref{small 4 chain}, $S = \{p_1, p_2, p_3, p_5\}$. The integer $n' = -1$ is an admissible number for $(Z, \omega, \pi, \theta)$ and $S$ because $S \subseteq \pi^{-1}([-1,0])$, and, as described in Example \ref{small 4 chain}, $\pi^{-1}(\{0\}) = \{p_1, p_2, p_3, p_4 \}$ and $\pi^{-1}(\{-1\}) = \{p_{-3}, p_{-2}, p_{-1}, p_0 \}$. Let $P' = \pi^{-1}([-1,0])$ and assume that $T_0$ is the following element of $\Tb(P', \omega)$.
\begin{center}
	$\tableau{8 & 7 & 6 & 4 \\ & 5 & 3 & 2 & 1}$
\end{center}
We illustrate $X_n(T_0)$, $X_{n,i}(T_0)$, and $v_n$ when $n = 1$. As described in Example \ref{matrix example}, $\Tb(S, \omega)$ has two elements. Hence, by Definition \ref{tableaux}, $|S, \omega| = 2$ and it follows that the range for the index $i$ is $1 \leq i \leq 2$. The $1^{st}$ set $X_1(T_0)$ derived from $(Z, \omega, \pi, \theta)$, $S$, and $T_0$ can be determined as follows. \\

Let $Q_1 = \pi^{-1}([-1,1])$ and let $P_1 = \pi^{-1}(\{1\})$. Consider the element $T \in \Tb(Q_1, \omega)$ depicted by the left-most diagram below. Using Lemma \ref{bipartition lemma} and setting $P = Q_1$, $S_0 = P'$, $S_1 = P \backslash P'$, and $\omega_P = \omega|_P$, it can be checked that the $1^{st}$ set $X_1(T_0)$ derived from $(Z, \omega, \pi, \theta)$, $S$, and $T_0$ is $X_1(T_0) = \{T\}$.
\begin{center}
	$\tableau{12 & 11 & 10 & 8 \\ & 9 & 7 & 6 & 5 \\ && 4 & 3 & 2 & 1}$ \;\;\;\;
	$\tableau{ {} & {} & {} & {} \\ & {} & {} & {} & \bullet \\ && {} & \bullet & \bullet & \bullet }$
\end{center}

Moreover, $\theta(S)$ is depicted by the cells in the right-most diagram above that are filled with bullets. Hence, with $T \in X_1(T_0)$ as above, $T|_{\theta(S)} \equiv_\theta R(1)$ where $R(1)$ is as described in Example \ref{matrix example}. So the $1^{st}$ part, $X_{1,1}(T_0)$, of the $1^{st}$ set derived from $(Z, \omega, \pi, \theta)$, $S$, $L$, and $T_0$, is $X_1(T_0)$. Moreover, the $2^{nd}$ part, $X_{1,2}(T_0)$, of the $1^{st}$ set derived from $(Z, \omega, \pi, \theta)$, $S$, $L$, and $T_0$, is the empty set because the only element of $X_1(T_0)$ is $T$. \\

Therefore, as $|X_{1,1}(T_0)| = 1$ and $|X_{1,2}(T_0)| = 0$, the $1^{st}$ vector derived from $(Z, \omega, \pi, \theta)$, $S$, $L$, and $T_0$ is the following column vector.
\begin{center}
	$v_1 = \begin{bmatrix}
	1 \\
	0
	\end{bmatrix}$
\end{center}
\end{example}

The sum of the entries in the column vectors defined in Definition \ref{combinatorial vectors} gives the number of periodic $P$-partitions.

\begin{lemma}\label{counting periodic shapes} Let $(Z, \omega, \pi, \theta)$ be a periodic quadruple system. Moreover, let $S$ be an index shape of $(Z, \omega, \pi, \theta)$, let $L$ be an indexing of $\Tb(S, \omega)$, let $R = L^{-1}$ be the inverse of $L$, let $P_n = \pi^{-1}([n])$ for all integers $n \geq 1$, let $n'$ be an admissible number for $(Z, \omega, \pi, \theta)$ and $S$, let $T_0 \in \Tb( \pi^{-1}([n',0]), \omega )$, and let $N = |S, \omega|$. Lastly, let $v_n$ be the $n^{th}$ vector derived from $(Z, \omega, \pi, \theta)$, $S$, $L$, and $T_0$. Then for all $n \geq 1$,
	$$|P_n, \omega| = \sum_{i=1}^N v_n(i). $$
\end{lemma}

\begin{example} Let $(Z, \omega, \pi, \theta)$ and $v_1$ be as described in Example \ref{combinatorial vector example}. The sum of the entries of $v_1$ is $1$. Moreover, as there is exactly one standard Young tableau on a single row, $|P_1, \omega| = 1$, where $P_1 = \pi^{-1}(\{1\})$. Hence, $|P_1, \omega|$ is also the sum of the entries of $v_1$.
\end{example}

\begin{proof} Let $n \in \mathbb{N}$, let $X_{n,i}(T_0)$ be the $i^{th}$ part of the $n^{th}$ set derived from $(Z, \omega, \pi, \theta)$, $S$, $L$, and $T_0$ for all $1 \leq i \leq N$, and let $X_n(T_0)$ be the $n^{th}$ set derived from $(Z, \omega, \pi, \theta)$, $S$, and $T_0$. Recall that, by Definition \ref{combinatorial vectors}, $X_{n,i}(T_0)$ is the set of elements $T \in X_n(T_0)$ such that
$$T|_{\theta^n(S)} \equiv_\theta R(i) $$
for all $1 \leq i \leq N$. Moreover, for all $T \in X_n(T_0)$, there is exactly one index $1 \leq i \leq N$ such that $T|_{\theta^n(S)} \equiv_\theta R(i)$ because $R(i_1) \neq R(i_2)$ for all $1 \leq i_1 \leq N$ and $1 \leq i_2 \leq N$ satisfying $i_1 \neq i_2$. Hence, we have
\begin{equation}\label{partitioning periodic p partitions}
X_n(T_0) = \bigcup_{i=1}^N X_{n,i}(T_0).
\end{equation}
Define the map $f : X_n(T_0) \to \Tb(P_n, \omega)$ by $f(T) = T|_{P_n}$. Let $P' = \pi^{-1}([n',0])$. For all $T \in X_n(T_0)$, for all $U \in T$, for all $p \in P'$, and for all $q \in P_n$, $U(p) > U(q)$. So if $T_1, T_2 \in X_n(T_0)$, then $T_1 \neq T_2$ if and only if $T_1|_{P_n} \neq T_2|_{P_n}$. It follows that $f$ is injective. Let $Q_n = P' \cup P_n$. Because $\pi : Z \to \mathbb{Z}$ is order-preserving, it follows that $P'$ is an order ideal of $Q_n$ and that $P_n$ is an order filter of $Q_n$. So by Lemma \ref{bipartition lemma} applied to $P = Q_n$, $S_0 = P'$, $S_1 = P_n$, and $\omega_P = \omega|_P$, $f$ is surjective. Hence, $f$ is a bijection, and the lemma follows from Equation \ref{partitioning periodic p partitions}.
\end{proof}

\subsection{The Matrix Difference Equation}

In this subsection, we prove that periodic $(P, \omega)$-partitions satisfy first order homogeneous matrix difference equations in which the matrices are the tableau transfer matrices. Then we prove that the periodic $(P, \omega)$-partitions can be enumerated with constant coefficient linear recurrence relations. As a consequence, we give a new proof of L\'opez et.al's and Sun's constant coefficient linear recurrence relation results by defining an appropriate periodic quadruple system $(Z, \omega, \pi, \theta)$, then explain how to prove a semi-standard variant of their results by modifying the labelling $\omega$ of $Z$ in our new proof. \\

We first prove a structural property of connected triples. It is a crucial lemma that is essential for what follows. Informally, the following lemma states that for connected triples of finite posets, we can define a notion of union for two $P$-partitions in a well-defined manner. \\

\begin{lemma}\label{extension} Let $Q$ be a finite poset, let $(A,B,C)$ be a connected triple of $Q$, and let $\omega_Q$ be a labelling of $Q$. Then for all $T' \in \Tb(A \cup B, \omega_Q)$ and $T'' \in \Tb(B \cup C, \omega_Q)$ such that $T'|_B = T''|_B$, there is a unique element $T \in \Tb(Q, \omega_Q)$ such that $T|_{A \cup B} = T'$ and $T|_{B \cup C} = T''$. 
\end{lemma}

\begin{example} Let $Q$ be the twelve element poset depicted by the below Young diagram. Moreover, let $A$ be the subposet of $Q$ depicted by the cells that are filled with asterisks, let $B$ be the subposet of $Q$ depicted by the cells that are filled with bullets, and let $C$ be the subposet of $Q$ depicted by the blank cells. Moreover, let $\omega_Q$ be a dual-natural labelling of $Q$, where in terms of tableaux the entries in the rows decreases from left to right, and the entries in the columns decrease from top to bottom.
\begin{center} 
	$\tableau{* & * & * & \bullet \\ & * & \bullet & \bullet & \bullet \\ && {} & {} & {} & {} }$
\end{center}
From how $A$, $B$, and $C$ are defined, $(A, B, C)$ is a connected triple of $Q$ by Definition \ref{connected triple}. Let $T'$ be the element of $\Tb(A \cup B, \omega_Q)$ that is depicted by the left-most diagram below, and let $T''$ be the element of $\Tb(B \cup C, \omega_Q)$ that is depicted by the right-most diagram below.
\begin{center}
	$\tableau{8 & 7 & 6 & 3 \\ & 5 & 4 & 2 & 1 } \;\;\;$
	$\tableau{ &&& 7 \\ && 8 & 5 & 3 \\ && 6 & 4 & 2 & 1 } \;\;\;$
\end{center}
From the above, we see that $T'|_B = T''_B$. Moreover, the unique element $T \in \Tb(Q, \omega_Q)$ such that  $T|_{A \cup B} = T'$ and $T|_{B \cup C} = T''$ is depicted by the following diagram.
\begin{center}
	$\tableau{12 & 11 & 10 & 7 \\ & 9 & 8 & 5 & 3 \\ && 6 & 4 & 2 & 1 }$
\end{center}
\end{example}

\begin{proof} It is enough to prove the following. Assume that there are elements $U_{1,1}, U_{1,2} \in \mathcal{A}(A \cup B, \omega_Q)$ and elements $U_{2,1}, U_{2,2} \in \mathcal{A}(B \cup C, \omega_Q)$ such that $U_{1,1} \equiv U_{1,2}$, $U_{2,1} \equiv U_{2,2}$, $U_{1,1}|_B \equiv U_{2,1}|_B$, and $U_{1,2}|_B \equiv U_{2,2}|_B$. Then the following two statements hold. \newline
	
\begin{enumerate}
	\item There exist elements $U_1 \in \mathcal{A}(Q, \omega_Q)$ and $U_2 \in \mathcal{A}(Q, \omega_Q)$ such that $U_1|_{A \cup B} \equiv U_{1,1}$, $U_1|_{B \cup C} \equiv U_{2,1}$, $U_2|_{A \cup B} \equiv U_{1,2}$, and $U_2|_{B \cup C} \equiv U_{2,2}$. \newline
	
	\item If $U_1' \in \mathcal{A}(Q, \omega_Q)$ and $U_2' \in \mathcal{A}(Q, \omega_Q)$ satisfy $U_1'|_{A \cup B} \equiv U_{1,1}$, $U_1'|_{B \cup C} \equiv U_{2,1}$, $U_2'|_{A \cup B} \equiv U_{1,2}$, and $U_2'|_{B \cup C} \equiv U_{2,2}$, then $U_1' \equiv U_2'$. \newline
\end{enumerate}

We first prove Statement 1. Since $U_{1,j}|_B \equiv U_{2,j}|_B$ for all $1 \leq j \leq 2$, there are order-embeddings $g_{1,j} : U_{1,j}(A \cup B) \to \mathbb{N}_0$ and $g_{2,j} : U_{2,j}(B \cup C) \to \mathbb{N}_0$ such that, for all $p \in B$,
$$g_{1,j}(U_{1,j}(p)) = g_{2,j}(U_{2,j}(p)).$$

So for all $1 \leq j \leq 2$, define $U_j : Q \to \mathbb{N}_0$ by
$$U_j(p) = \begin{cases}
g_{1,j}(U_{1,j}(p)) & \text{if }p \in A \cup B \\
g_{2,j}(U_{2,j}(p)) & \text{if }p \in B \cup C. \\
\end{cases} $$

For all $1 \leq j \leq 2$, the above map $U_j$ is well-defined because of the definition of $g_{1,j}$ and $g_{2,j}$. Moreover, for all $1 \leq j \leq 2$, the map $U_j$ satisfies 
\begin{center}
	$U_j|_{A \cup B} = g_{1,j} \circ U_{1,j} \equiv U_{1,j}$ and $U_j|_{B \cup C} = g_{2,j} \circ U_{2,j} \equiv U_{2,j}$.
\end{center}
Hence, to prove Statement 1, it is enough to prove that $U_j \in \mathcal{A}(Q, \omega_Q)$ for all $1 \leq j \leq 2$. So let $j \in \{1,2\}$. To see that $U_j$ is order reversing as required by the definition of $\mathcal{A}(Q, \omega_Q)$ in Section \ref{sec:prelims}, suppose otherwise. \\

Because $U_j|_{A \cup B} \in \mathcal{A}(A \cup B, \omega_Q)$ and $U_j|_{B \cup C} \in \mathcal{A}(B \cup C, \omega_Q)$, $U_j|_{A \cup B}$ and $U_j|_{B \cup C}$ are order reversing. \\

So, as we are supposing that $U_j$ is not order reversing, there are elements $p \in A$ and $q \in C$ such that $p < q$ but $U_j(p) < U_j(q)$. By Property 2 of Definition \ref{connected triple}, there exists an element $p' \in B$ such that $p < p' < q$ in $Q$. As $p, p' \in A \cup B$, as $p < p'$ in $A \cup B$, and as $U_j|_{A \cup B}$ is order reversing, it follows that $U_j(p) \geq U_j(p')$. Moreover, as $p', q \in B \cup C$, as $p' < q$ in $B \cup C$, and as $U_j|_{B \cup C}$ is order reversing, it follows that $U_j(p') \geq U_j(q)$. But then, $U_j(p) \geq U_j(p') \geq U_j(q)$, which is contrary to the assumption that $U_j(p) < U_j(q)$. \\

So $U_j$ is order reversing. Suppose that $U_j \notin \mathcal{A}(Q, \omega_Q)$. Then there are elements $p, q \in Q$ such that $p < q$ in $Q$, $\omega_Q(p) > \omega_Q(q)$, and $U_j(p) = U_j(q)$. Because $U_j|_{A \cup B} \in \mathcal{A}(A \cup B, \omega_Q)$ and $U_j|_{B \cup C} \in \mathcal{A}(B \cup C, \omega_Q)$, it follows that $p \in A$ and $q \in C$. Since $p \in A$ and $q \in C$, Property 2 of Definition \ref{connected triple} implies that there exists an element $p' \in B$ such that $p < p' < q$ in $Q$. \\

If $\omega_Q(p) \leq \omega_Q(p')$ and $\omega_Q(p') \leq \omega_Q(q)$, then $\omega_Q(p) \leq \omega_Q(p') \leq \omega_Q(q)$, implying that $\omega_Q(p) \leq \omega_Q(q)$. But that is contrary to the assumption that $\omega_Q(p) > \omega_Q(q)$. Hence, $\omega_Q(p) > \omega_Q(p')$ or $\omega_Q(p') > \omega_Q(q)$. So as $p,p' \in A \cup B$, $p < p'$ in $A \cup B$, $p',q \in B \cup C$, $p' < q$ in $B \cup C$, $U_j|_{A \cup B} \in \mathcal{A}(A \cup B, \omega_Q)$, and $U_j|_{B \cup C} \in \mathcal{A}(B \cup C, \omega_Q)$, it follows that 
\begin{center}
	$U_j(p) > U_j(p') \geq U_j(q)$ \; or \; $U_j(p) \geq U_j(p') > U_j(q)$.
\end{center}
But then, $U_j(p) > U_j(q)$, which is contrary to the assumption that $U_j(p) = U_j(q)$. Hence, $U_j \in \mathcal{A}(Q, \omega_Q)$ and Statement 1 follows. \\

To prove Statement 2, let $U_{1,1}$, $U_{1,2}$, $U_1'$, $U_{1,2}$, $U_{2,2}$, and $U_2'$ be as described in the beginning of the proof, and suppose that $U_1'$ is not order equivalent to $U_2'$. Because $U'_1|_{A \cup B} \equiv U_{1,1} \equiv U_{1,2} \equiv U_2'|_{A \cup B}$ and because $U_1'|_{B \cup C} \equiv U_{2,1} \equiv U_{2,2} \equiv U_2'|_{B \cup C}$, we have that $U_1'|_{A \cup B} \equiv U_2'|_{A \cup B}$ and $U_1'|_{B \cup C} \equiv U_2'|_{B \cup C}$. So there are elements $p,q \in Q$ such that $p \in A$, $q \in C$, and exactly one of the following holds. \newline

\begin{itemize}
	\item $U_1'(p) < U_1'(q)$ and $U_2'(p) > U_2'(q)$ \newline
	\item $U_1'(p) > U_1'(q)$ and $U_2'(p) < U_2'(q)$ \newline
	\item $U_1'(p) = U_1'(q)$ and $U_2'(p) \neq U_2'(q)$ \newline
	\item $U_1'(p) \neq U_1'(q)$ and $U_2'(p) = U_2'(p)$ \newline
\end{itemize}

Suppose that $U_1'(p) < U_1'(q)$ and $U_2'(p) > U_2'(q)$. Then, as $U_1'$ and $U_2'$ are order reversing maps, it follows that $p \parallel q$ in $Q$. But as $(A, B, C)$ is a connected triple of $Q$, that violates Property 2 of Definition \ref{connected triple}. Similarly, if $U_1'(p) > U_1'(q)$ and $U_2'(p) < U_2'(q)$, then Property 2 of Definition \ref{connected triple} would be violated. So without loss of generality, suppose that $U_1'(p) = U_1'(q)$ and $U_2'(p) \neq U_2'(q)$. \\

Since $p \in A$ and $q \in C$, Property 2 of Definition \ref{connected triple} implies that there exists an element $p' \in B$ such that $p < p' < q$ in $Q$. Hence, as $U_1'$ is order reversing and as $U_1'(p) = U_1'(q)$, $U_1'(p) = U_1'(p') = U_1'(q)$. So as $p, p' \in A \cup B$ and $U_1'|_{A \cup B} \equiv U_2'|_{A \cup B}$, it follows that $U_2'(p) = U_2'(p')$. Similarly, as $p', q \in B \cup C$ and $U_1'|_{B \cup C} \equiv U_2'|_{B \cup C}$, it follows that $U_2'(p') = U_2'(q)$. But then, $U_2'(p) = U_2'(q)$, which is contrary to the assumption that $U_2'(p) \neq U_2'(q)$. Hence, Statement 2 follows.
	
\end{proof}

\begin{remark}\label{info} Note that the converse of Lemma \ref{extension} is also true. If $T \in \Tb(Q, \omega)$, then $T$ uniquely determines $T|_{A \cup B} \in \Tb(A \cup B, \omega)$ and $T|_{B \cup C} \in \Tb(B \cup C, \omega)$. 
\end{remark}

In preparation for the main result of this section, we introduce the following technical definition. It defines a positive integer that depends on the periodic quadruple system being considered.

\begin{definition}\label{minimum number} Let $(Z, \omega, \pi, \theta)$ be a periodic quadruple system and let $S$ be an index shape of $(Z, \omega, \pi, \theta)$. Then the \emph{minimum number of $(Z, \omega, \pi, \theta)$ and $S$} is the smallest positive integer $m$ such that if $n \in \mathbb{N}$ satisfies $n \geq m$, then for all $p \in Z$ satisfying $\pi(p) = n+1$, there exists an element $q \in \theta^{n}(S)$ such that $q < p$ in $Z$ and $\pi(q) \geq 1$.
\end{definition}

\begin{example}\label{minimum ex} Let $(Z, \omega, \pi, \theta)$ and $S$ be as in Example \ref{small 4 chain}. We explain why the minimum number $m$ for $(Z, \omega, \pi, \theta)$ and $S$ is $1$. Let the twelve-element poset $\pi^{-1}([0,2])$ be depicted by the Young diagram below. Moreover, let $\pi^{-1}(\{0\})$ be depicted by the four cells in the top row of the diagram, let $\theta(S)$ be depicted by the cells filled with asterisks or circles, and let $\pi^{-1}(\{2\})$ be depicted by the cells filled with bullets.
\begin{center}
	$\tableau{{} & {} & {} & * \\ & {} & \circ & \circ & \circ \\ && \bullet & \bullet & \bullet & \bullet}$
\end{center}
Set $m = 1$ and set $n = m$. From the above diagram, it can be checked that if $p \in Z$ satisfies $\pi(p) = n + 1 = 2$, then $p$ is represented by one of the cells filled with bullets. Fixing any such element $p \in Z$ satisfying $\pi(p) = 2$, it can also be checked, from the same diagram, that there exists an element $q \in \theta^{n}(S) = \theta(S)$, specifically one of the cells filled with a circle, such that $q < p$ in $Z$ and $\pi(q) \geq 1$. Hence, as it can be checked that the above is also true when $m = 1$ and $n > 1$, it follows that $m = 1$ is the minimum number of $(Z, \omega, \pi, \theta)$.
\end{example}

\begin{remark}\label{admissible bound} The minimum number of a periodic quadruple system always exists. Let $(Z, \omega, \pi, \theta)$ be a periodic quadruple system and let $S$ be an index shape of $(Z, \omega, \pi, \theta)$. Since $S$ is finite, there is a positive integer $n$ such that, for all $p \in \theta^n(S)$, $\pi(p) \geq 1$. Moreover, $S$ connects $Z$, so as $\theta$ is an order-automorphism on $Z$, $\theta^n(S)$ connects $Z$. Hence, by Property 2 of Definition \ref{connected triple}, it follows that for all $p \in Z$ satisfying $\pi(p) = n+1$ and for all $p' \in Z$ satisfying $\pi(p') \leq 0$, there exists an element $q \in \theta^n(S)$ such that $p' < q < p$. But as $\pi(q) \geq 1$ for all $q \in \theta^n(S)$, it follows from Definition \ref{minimum number} that the minimum number of $(Z, \omega, \pi, \theta)$ exists and is at most $n$.
\end{remark}

We now prove that periodic $P$-partitions satisfy a first order homogeneous matrix difference equation involving tableau transfer matrices. For the proof, recall the definitions of the sets $X_n(T_0)$ and $X_{n,i}(T_0)$ in Definition \ref{combinatorial vectors}.

\begin{theorem}\label{difference equation} Let $(Z, \omega, \pi, \theta)$ be a periodic quadruple system, let $S$ be an index shape of $(Z, \omega, \pi, \theta)$, and let $L$ be an indexing of $\Tb(S, \omega)$. Moreover, let $n'$ be an admissible number for $(Z, \omega, \pi, \theta)$ and $S$, let $T_0 \in \Tb( \pi^{-1}([n',0]), \omega)$, and let $(v_n)_{n=1,2,\dots}$ be a sequence such that for all $n \geq 1$, $v_n$ is the $n^{th}$ vector derived from $(Z, \omega, \pi, \theta)$, $S$, $L$, and $T_0$. Lastly, let $M$ be the tableau transfer matrix derived from $(Z, \omega, \pi, \theta)$, $S$, and $L$, and let $m$ be the minimum number of $(Z, \omega, \pi, \theta)$ and $S$. Then for all $n \geq m$,
$$v_{n+1} = M \, v_n. $$
\end{theorem}

\begin{example}\label{difference eq ex} Let $(Z, \omega, \pi, \theta)$, $S$, $L$, $T_0$ be as in Example \ref{combinatorial vector example}. By Example \ref{minimum ex}, the minimum number $m$ for $(Z, \omega, \pi, \theta)$ and $S$ is $1$. So consider the $1^{st}$ vector $v_1$ derived from $(Z, \omega, \pi, \theta)$, $S$, and $L$. As shown in Example \ref{combinatorial vector example},
$$v_1 = \begin{bmatrix} 1 \\ 0 \end{bmatrix}.$$
Next, let $M$ be the tableau transfer matrix from Example \ref{matrix example}. 
What this theorem allows us to do is to determine $v_n$ from $M$ and $v_1$ for any $n \geq 1$. For instance,
$$ v_2 = M v_1 = \begin{bmatrix}
3 & 4 \\
2 & 3 \\
\end{bmatrix} \; \begin{bmatrix} 1 \\ 0 \end{bmatrix} = \begin{bmatrix} 3 \\ 2 \end{bmatrix} . $$
\end{example}

\begin{proof} Let $n \in \mathbb{N}$. Define $Q = \pi^{-1}([n',n+1])$, $C = \pi^{-1}(\{n+1\})$, $B = \theta^n(S)$, and $A = \pi^{-1}([n',n]) \backslash B$. Then $(A, B, C)$ is a connected triple of $Q$. Note that $A \cup B = \pi^{-1}([n',n])$ and $\theta(B) = \theta^{n+1}(S)$. \\

If $T \in X_{n+1}(T_0)$, then $T|_{A \cup B} \in \Tb(A \cup B, \omega)$ and $T|_{B \cup C} \in \Tb(B \cup C, \omega)$. Since $T \in X_{n+1}(T_0)$, Definition \ref{combinatorial vectors} implies that for all $U \in T$, $p \in \pi^{-1}([n',0])$, and $q \in \pi^{-1}([n+1])$, $U(p) > U(q)$. It follows that $T|_{A \cup B} \in X_n(T_0)$  and $T|_{B \cup C} \in \Tb(B \cup C, \omega)$. Moreover, see Remark \ref{info}, $T|_{A \cup B}$ and $T|_{B \cup C}$ are uniquely determined by $T \in X_{n+1}(T_0)$.  \\ 

Next, define the map 
$$f : X_{n+1}(T_0) \to \{(T'', T') \in \Tb(B \cup C, \omega) \times X_n(T_0) : T''|_B = T'|_B \}$$ 
by
$$f(T) = (T|_{B \cup C}, T|_{A \cup B})$$
for all $T \in X_{n+1}(T_0)$. By what we just showed, $f$ is well-defined and injective. We will prove that $f$ is also a bijection. To that end, it is enough to show that $f$ is surjective. \\

Let $T' \in X_n(T_0)$, let $T'' \in \Tb(B \cup C, \omega)$, and assume that $T'|_B = T''|_B$. Recall that $\Tb(Q, \omega) = \Tb(Q, \omega_Q)$, where $\omega_Q$ is the labelling of $Q$ such that $\omega|_Q \equiv \omega_Q$. Hence, by Lemma \ref{extension}, there is a unique element $T \in \Tb(Q, \omega)$ such that $T|_{A \cup B} = T'$ and $T|_{B \cup C} = T''$. So to prove that $f$ is surjective, it is enough to prove that $T \in X_{n+1}(T_0)$. \\

Let $P' = \pi^{-1}([n',0])$ and let $P = \pi^{-1}([n+1])$. Because $T|_{A \cup B} = T' \in X_n(T_0)$ and because $P' \subseteq A \cup B$, Definition \ref{combinatorial vectors} implies that $T|_{P'} = T_0$. Hence, by Definition \ref{combinatorial vectors}, it is enough to show that for all $U \in T$, $p \in P'$, and $q \in P$, $U(p) > U(q)$. To that end, let $U \in T$, let $p \in P'$, and let $q \in P$. \\

If $\pi(q) \leq n$, then $p,q \in A \cup B$, implying, as $T|_{A \cup B} \in X_n(T_0)$, that $U(p) > U(q)$ by Definition \ref{combinatorial vectors} applied to $T|_{A \cup B}$. So assume without loss of generality that $\pi(q) = n+1$. Because $n \geq m$, by hypothesis where $m$ is the minimum number of $(Z, \omega, \pi, \theta)$ and $S$, Definition \ref{minimum number} implies that there exists an element $p' \in \theta^n(S)$ such that $p' < q$ in $Z$ and $\pi(p') \geq 1$. In particular, as $p' < q$ and $U$ is order reversing, $U(p') \geq U(q)$. \\

Because $B = \theta^n(S)$, $p' \in B$. Moreover, $T|_{A \cup B} \in X_n(T_0)$ and $p \in P'$. Furthermore, by Definition \ref{combinatorial vectors}, $U(p'') > U(q'')$ for all $p'' \in P'$ and $q'' \in \pi^{-1}([n])$. Lastly, $p' \in \pi^{-1}([n])$ because $p' \in \theta^n(S)$ and $\pi(p') \geq 1$. So as $p \in P'$ and $p' \in \pi^{-1}([n])$, we have $U(p) > U(p')$. Hence,
$$U(p) > U(p') \geq U(q),$$
implying that $U(p) > U(q)$. From this, it follows that $T \in X_{n+1}(T_0)$. Hence, $f$ is a bijection. \\

Whence, for all $T_2 \in \Tb(S, \omega)$, the number of elements $T \in X_{n+1}(T_0)$ satisfying $T|_{\theta^{n+1}(S)} \equiv_\theta T_2$ is
\begin{equation}\label{intermediary step}
\sum_{T_1 \in \Tb(S, \omega)}  M(T_2, T_1)\; |\{T \in X_n(T_0) : T|_{\theta^n(S)} \equiv_\theta T_1 \}|
\end{equation}
for the following reasons. \\

Since $B = \theta^n(S)$, the fact that $f$ is a bijection implies that the following is true for all $T_2 \in \Tb(S, \omega)$. If $g$ is the restriction of $f$ to the set of elements $T \in X_{n+1}(T_0)$ satisfying $T|_{\theta^{n+1}(S)} \equiv_\theta T_2$, then $g$ is injective and the range of $g$ is the set of ordered pairs $(T'', T') \in \Tb(B \cup C, \omega) \times X_n(T_0)$ satisfying $T''|_{\theta^{n+1}(S)} \equiv_\theta T_2$ and $T''|_{\theta^n(S)} = T'|_{\theta^n(S)}$. \\

Moreover, as $\theta$ is an order-automorphism on $Z$, as 
$$B \cup C = \theta^n(S) \cup \pi^{-1}(\{n+1\}) = \theta^{n+1}(\theta^{-1}(S) \cup \pi^{-1}(\{0\})),$$
as $\theta^n(S) = \theta^{n+1}(\theta^{-1}(S))$, and as $M$ is the tableau transfer matrix derived from $(Z, \omega, \pi, \theta)$, $S$ and $L$, Definition \ref{tableau transfer matrix} implies that for all $T_1,T_2 \in \Tb(S, \omega)$, the number of elements $T'' \in \Tb(B \cup C, \omega)$ satisfying $T''|_{\theta^{n+1}(S)} \; \equiv_\theta \; T_2$ and $T''|_{\theta^n(S)} \equiv_\theta T_1$ is $M(T_2, T_1)$. \\

Hence, the number of elements $T \in X_{n+1}(T_0)$ satisfying $T|_{\theta^{n+1}(S)} \equiv_\theta T_2$ is given by Expression \ref{intermediary step}. \\

Lastly, by Definition \ref{combinatorial vectors},
$$ |X_{n, L(T_1)}(T_0)| = |\{T \in X_n(T_0) : T|_{\theta^n(S)} \equiv_\theta T_1 \}| $$
for all $T_1 \in \Tb(S, \omega)$, and
$$ |X_{n+1, L(T_2)}(T_0)| = |\{T \in X_{n+1}(T_0) : T|_{\theta^{n+1}(S)} \equiv_\theta T_2 \}| $$
for all $T_2 \in \Tb(S, \omega)$. \\

Therefore, for all $T_2 \in \Tb(S, \omega)$, Definition \ref{combinatorial vectors} implies that
\begin{align*}
v_{n+1}(L(T_2)) &= |X_{n+1, L(T_2)}(T_0)| \\
&= |\{T \in X_{n+1}(T_0) : T|_{\theta^{n+1}(S)} \equiv_\theta T_2\}| \\
&= \sum_{T_1 \in \Tb(S, \omega)}  M(T_2, T_1)\; |\{T \in X_n(T_0) : T|_{\theta^n(S)} \equiv_\theta T_1 \}| \\
&= \sum_{T_1 \in \Tb(S, \omega)}  M(T_2, T_1)\; |X_{n, L(T_1)}(T_0)| \\
&= \sum_{T_1 \in \Tb(S, \omega)} M(T_2,T_1) \; v_n(L(T_1)).
\end{align*}
From this, the theorem follows from the definition of matrix multiplication.
\end{proof}

We now establish the following enumeration result for periodic $P$-partitions.

\begin{theorem}\label{recurrence relations} Let $(Z, \omega, \pi, \theta)$ be a periodic quadruple system, let $S$ be an index shape of $(Z, \omega, \pi, \theta)$, let $L$ be an indexing of $\Tb(S, \omega)$, let $M$ be the tableau transfer matrix derived from $(Z, \omega, \pi, \theta)$, $S$, and $L$, let $m$ be the minimum number of $(Z, \omega, \pi, \theta)$ and $S$, let $\nu(x)$ be the minimal polynomial of $M$, and let $d$ be the largest non-negative integer $d$ such that $x^d$ divides $\nu(x)$. Lastly, let $P_n = \pi^{-1}([n])$ for all integers $n \geq 1$. Then the sequence $(a_n)_{n=m+d,m+d+1,\dots}$ defined by $a_n = |P_n, \omega|$ for all $n \geq m+d$, satisfies a constant coefficient linear recurrence whose characteristic polynomial is $\nu(x) / x^d$.
\end{theorem}

\begin{example}\label{recurrence relations ex} Let $(Z, \omega, \pi, \theta)$, $S$, $L$, $T_0$, and $M$ be as in Example \ref{difference eq ex}. Moreover, let $(a_n)_{n=1,2,\dots}$ be defined by $a_n = |P_n, \omega|$ for all $n \geq 1$ where $P_n = \pi^{-1}([n])$. The minimal polynomial $\nu(x)$ of $M$ is $x^2 - 6x + 1$; in particular, $d = 0$ and $\nu(x)/x^d = x^2 - 6x + 1$. Moreover, by Example \ref{minimum ex}, the minimum number $m$ of $(Z, \omega, \pi, \theta)$ and $S$ is $1$. So, as $d = 0$ and $m = 1$, the sequence $(a_n)_{n=1,2,\dots}$ satisfies the recurrence
$$a_{n+2} = 6 a_{n+1} - a_n. $$
This gives a new proof of one of the recurrence relation results from Sun \cite{EYTPS} and L\'opez et.al \cite{CHC}.
\end{example}

\begin{proof} Let $N = |S, \omega|$, let $n'$ be an admissible number for $(Z, \omega, \pi, \theta)$ and $S$, and let $T_0 \in \Tb(\pi^{-1}([n',0])$. Since $\nu(x)$ is the minimal polynomial of $M$, $\nu(M) = 0$. So by Theorem \ref{difference equation},
 $$\nu(M) \, v_n = 0$$
for all $n \geq m$.	Hence, the theorem follows from Lemma \ref{counting periodic shapes}.
\end{proof}

We now give new proofs of the constant coefficient linear recurrence relations from \cite{CHC, EYTPS, HAlgebra}. For all integers $n \geq 1$ and $k \geq 1$, let $X_{n,k}$ denote the parallelogramic shape with $n$ rows and $k$ cells in each row. Moreover, let $T(n,k)$ denote the number of standard Young tableaux of shape $X_{n,k}$.

\begin{corollary}\label{existence}[L\'opez et.al. (2017, \cite{CHC}), cf. Sun \cite{EYTPS}, \cite{HAlgebra}] Fix a positive integer $k \geq 3$. Moreover, let $a_n = T(n,k)$ for all $n \geq 1$. Then the sequence $(a_n)_{n=1,2,\dots}$ satisfies a constant coefficient linear recurrence.
\end{corollary}

\begin{figure}
\begin{center}
$\tableau{& {} \\ & {} & {} \\ {} & {} & {} & {}}$ \;\;\;\;\;\;\;\;\;\;
$\begin{bmatrix}
4 & 5 & 5 & 6 & 6 & 7 & 7 \\
3 & 4 & 4 & 5 & 5 & 6 & 6 \\
3 & 4 & 4 & 5 & 5 & 6 & 6 \\
2 & 3 & 3 & 4 & 4 & 5 & 5 \\
2 & 3 & 3 & 4 & 4 & 5 & 5 \\
0 & 0 & 2 & 0 & 3 & 0 & 4 \\
0 & 0 & 2 & 0 & 3 & 0 & 4 
\end{bmatrix}$
\end{center}
\caption{}\label{fig:five-shape}
\end{figure}

\begin{example}\label{k eq 5 example} For $k = 5$, consider the periodic quadruple system $(Z(X, \Delta), \omega, \pi, \theta)$ from Proposition \ref{Full Euclidean system} when $X = \{(1,0), (2,0), (3,0), (4,0), (5,0) \}$ and $\Delta = (1,1)$. Moreover, let
$$S = \{(1,0), (2,0), (3,0), (4,0) \} \cup \{(2,-1), (3,-1) \} \cup \{(2,-2) \}. $$
This is depicted by the seven-cell Young diagram in Figure \ref{fig:five-shape}. It can be checked that $S$ is an index shape of $(Z(X, \Delta), \omega, \pi, \theta)$. Moreover, there is a labelling $L$ of $\Tb(S, \omega)$ such that the tableau transfer matrix derived from $(Z(X, \Delta), \omega, \pi, \theta)$, $S$, and $L$ is the matrix in Figure \ref{fig:five-shape}. The characteristic polynomial of this matrix is $x^4 - 24x^3 + 40x^2 + 8x$, so by the Cayley-Hamilton Theorem, the minimal polynomial of this matrix divides the aforementioned polynomial. Lastly, it can be checked, in a manner similar to what was done in Example \ref{minimum ex}, that the minimum number $m$ of $(Z, \omega, \pi, \theta)$ and $S$ is $m = 1$. Hence, by Theorem \ref{recurrence relations}, the subsequence $(a_n)_{n=2,3,\dots}$ of the sequence $(a_n)_{n=1,2,\dots}$ satisfies a constant coefficient linear recurrence. So from Hardin and Heinz's verification that $T(1,5)$, $T(2,5)$, $T(3,5)$, and $T(4,5)$ satisfy the above recurrence, we obtain a new proof of the recurrence relation empirically observed by Hardin and Heinz (\cite{OEIS}, A181196) for $T(n,5)$.
\end{example}

\begin{proof} Consider the periodic quadruple system $(Z(X, \Delta), \omega, \pi, \theta)$ from Proposition \ref{Full Euclidean system} when $X = \{(i,0) : 1 \leq i \leq k\}$ and $\Delta = (1,1)$. Next, let $S$ be defined by
$$ S = \{(i,0) : i \in [n] \} \cup  \bigcup_{k=1}^{n-2} \{(i, -k) : i \in [2, n-k] \}. $$
It can be checked that $S$ is an index shape of $Z(X, \Delta)$, so the corollary follows from Theorem \ref{recurrence relations}.
\end{proof}

Lastly, we observe that by replacing the labelling $\omega$ of $Z(X, \Delta)$ in the proof of Corollary \ref{existence} with a labelling $\omega''$ of $Z(X, \Delta)$ that is analogous to the labelling $\omega''$ in Example \ref{small 4 chain}, we also obtain the following corollary.

\begin{corollary}\label{semi-standard existence} Fix a positive integer $k \geq 3$. Moreover, let $b_n$ denote the number of semi-standard tableau classes on $X_{n,k}$ for all $n \geq 1$. Then the sequence $(b_n)_{n=1,2,\dots}$ satisfies a constant coefficient linear recurrence relation.
\end{corollary}

\section{Asymptotic Growth Rates} \label{sec: asymptotic}

In this section, we prove that many tableau transfer matrices are primitive. Using this fact, we prove that the asymptotic growth rate of the number $|P, \omega|$ of periodic $(P, \omega)$-partitions is given by $c \, r^n$, where $c$ and $r$ are constants. Afterwards, we briefly remark on how our proof that the above matrices are primitive can give further descriptions for the constant coefficient linear recurrence relations we are analysing. \\

Frobenius \cite{MA} introduced the notion of a \emph{primitive matrix} for Perron-Frobenius theory and established that such matrices are equivalent to the following, which we give as a definition.

\begin{definition}[Frobenius \cite{MA}] A \emph{primitive matrix} is a square matrix $M$ such that $M$ is a non-negative matrix and such that, for some integer $k \geq 1$, $M^k$ is a positive matrix.
\end{definition}

Motivated by the above definition, we define the following index shapes and tableau transfer matrices.

\begin{definition}\label{def:separable index shape} Let $(Z, \omega, \pi, \theta)$ be a periodic quadruple system and let $S$ be an index shape of $(Z, \omega, \pi, \theta)$. Then $S$ is a \emph{separable index shape of $(Z, \omega, \pi, \theta)$} if $S$ is an order filter of $\pi^{-1}(\mathbb{Z} \backslash \mathbb{N})$ and if $S$ is an order ideal of $S \cup \pi^{-1}(\mathbb{N})$. When $(Z, \omega, \pi, \theta)$ is not specified, we call $S$ a \emph{separable index shape}.
\end{definition}

\begin{remark}\label{rem:separable index shapes always exist} The index shape constructed in Remark \ref{rem:index shapes always exist} is also a separable index shape. From this, it follows that separable index shapes always exist.
\end{remark}

\begin{example} The index shapes described in Example \ref{small 4 chain} and Example \ref{k eq 5 example}, and the index shapes described in the proof of Corollary \ref{existence} are examples of separable index shapes.
\end{example}

\begin{definition}\label{def:separable matrices} Let $(Z, \omega, \pi, \theta)$ be a periodic quadruple system, let $S$ be an index shape of $(Z, \omega, \pi, \theta)$, and let $M$ be a tableau transfer matrix derived from $(Z, \omega, \pi, \theta)$ and $S$. Then call $M$ \emph{separable} if $S$ is separable index shape of $(Z, \omega, \pi, \theta)$.
\end{definition}

\begin{example}\label{ex:separable matrices} The tableau transfer matrices in Example \ref{matrix example} and Example \ref{k eq 5 example}, and the tableau transfer matrices described in the proof of Corollary \ref{existence} are examples of tableau transfer matrices that are separable.
\end{example}

We prove the following property of separable tableau transfer matrices.

\begin{lemma}\label{tableau transfer matrices primitive} Any tableau transfer matrix that is separable is also primitive.
\end{lemma}

\begin{example} The matrices in Example \ref{ex:separable matrices} are tableau transfer matrices that are primitive.
\end{example}

\begin{proof} Let $(Z, \omega, \pi, \theta)$ be a periodic quadruple system, let $S$ be a separable index shape of $(Z, \omega, \pi, \theta)$, let $L$ be a labelling of $\Tb(S, \omega)$, let $M$ be the tableau transfer matrix derived from $(Z, \omega, \pi, \theta)$, $S$, and $L$, let $R = L^{-1}$ be the inverse of $L$, and let $N$ denote the number of rows of $M$. Because $S$ is a separable index shape of $(Z, \omega, \pi, \theta)$, $M$ is separable. By Definition \ref{tableau transfer matrix}, any tableau transfer matrix is a non-negative square matrix. Hence, it is enough to prove that $M^k$ is positive for some integer $k \geq 1$. We will first prove by induction on $k$ that for all $1 \leq i \leq N$ and $1 \leq j \leq N$, $M^k(i,j)$ equals to the number of elements $T \in \Tb(\theta^{-k}(S) \cup \pi^{-1}([-k+1,0]), \omega)$ such that $T|_S \equiv_\theta R(i)$ and $T|_{\theta^{-k}(S)} \equiv_\theta R(j)$. \\
	
The base case $k = 1$ follows from Definition \ref{tableau transfer matrix}. So let $k \geq 1$ and assume that the claim holds for all integers in $[k]$. Consider the matrix $M^{k+1}$. Then for all $1 \leq i \leq N$ and for all $1 \leq j \leq N$,
\begin{equation}\label{component relation}
M^{k+1}(i,j) = \sum_{i' = 1}^N M^k(i,i') \, M(i',j).
\end{equation}
For all indices $i'$, $i$, and $j$ in the summation in Equation \ref{component relation}, we note the following. Definition \ref{tableau transfer matrix}, and the fact that $\theta$ is an order-automorphism on $Z$, implies that $M(i',j)$ is the number of elements $T \in Tb(\theta^{-k-1}(S) \cup \pi^{-1}(\{-k\}), \omega)$ such that $T|_{\theta^{-k}(S)} \equiv_\theta R(i')$ and $T|_{\theta^{-k-1}(S)} \equiv_\theta R(j)$. By the induction hypothesis applied to the entries $M^k(i,i')$ of $M^k$, $M^{k}(i,i')$ equals to the number of elements $T \in \Tb(\theta^{-k}(S) \cup \pi^{-1}([-k+1,0]), \omega)$ satisfying $T|_S \equiv_\theta R(i)$ and $T|_{\theta^{-k}(S)} \equiv_\theta R(i')$. So as
$$((\theta^{-k-1}(S) \cup \pi^{-1}(\{-k\})) \backslash \theta^{-k}(S), \theta^{-k}(S), \pi^{-1}([-k+1,0]))$$
is a connected triple of $\theta^{-k-1}(S) \cup \pi^{-1}([-k,0])$, it follows from Lemma \ref{extension}, Remark \ref{info}, and Equation \ref{component relation} that $M^{k+1}(i,j)$ equals to the number of elements $T \in \Tb(\theta^{-k-1}(S) \cup \pi^{-1}([-k,0]), \omega)$ that satisfy $T|_S \equiv_\theta R(i)$ and $T|_{\theta^{-k-1}(S)} \equiv_\theta R(j)$. From this, the claim follows. Note that the claim just proved does not depend on the assumption that $S$ is a separable index shape of $(Z, \omega, \pi, \theta)$. \\

Next, we prove the following. By Definition \ref{periodic quadruple system}, $\pi(\theta(p)) = \pi(p) + 1$ for all $p \in Z$. So as $S$ and $\pi(S)$ are finite, there exists a sufficiently large integer $k \geq 1$ such that $\pi(S) \cap \pi(\theta^{-k}(S)) = \emptyset$. Hence, for some integer $k \geq 1$,
\begin{equation}\label{disjoint relation}
S \cap \theta^{-k}(S) = \emptyset.
\end{equation}
Let $S^* = \theta^{-k}(S) \cup \pi^{-1}([-k+1,0])$ and let $X_k = S^* \backslash (S \cup \theta^{-k}(S))$. By Property 2 of Definition \ref{index shapes} and the fact that $\theta$ is an order-automorphism on $Z$, $S \subseteq S^*$ since
$$S \subseteq \theta^{-1}(S) \cup \pi^{-1}(\{0\}) \subseteq \theta^{-2}(S) \cup \pi^{-1}(\{-1,0\}) \subseteq \cdots \subseteq \theta^{-k}(S) \cup \pi^{-1}([-k+1,0]) = S^*.$$
Because $S$ is a separable index shape of $(Z, \omega, \pi, \theta)$, Definition \ref{def:separable index shape} implies that $S$ is an order ideal of $S \cup \pi^{-1}(\mathbb{N})$. So as $\theta$ is an order-automorphism on $Z$, $\theta^{-k}(S)$ is an order ideal of $S^*$. Moreover, as $S$ is a separable index shape of $(Z, \omega, \pi, \theta)$, it follows from Definition \ref{def:separable index shape} that $S$ is an order filter of $\pi^{-1}(\mathbb{Z} \backslash \mathbb{N})$. This fact, combined with Equation \ref{disjoint relation}, implies that $S$ is an order filter of $X_k \cup S$. So as Equation \ref{disjoint relation} is true, the following holds. By applying Lemma \ref{bipartition lemma} to $S^*$, $\omega|_{S^*}$, $\theta^{-k}(S)$, and $X_k \cup S$, then by applying, if $X_k \neq \emptyset$, Lemma \ref{bipartition lemma} to $X_k \cup S$, $\omega|_{X_k \cup S}$, $X_k$, and $S$, we conclude that there exists an element $T \in \Tb(S^*, \omega)$ such that $T|_{S} \equiv_\theta R(i)$ and $T|_{\theta^{-k}(S)} \equiv_\theta R(j)$. Hence, $M(i,j) \geq 1$ and the lemma follows.
\end{proof}

A consequence of the above is that we can combinatorially describe the entries of any power of a tableau transfer matrix.

\begin{corollary}\label{powers} Let $(Z, \omega, \pi, \theta)$ be a periodic quadruple system, let $S$ be an index shape of $(Z, \omega, \pi, \theta)$, let $L$ be an indexing of $\Tb(S, \omega)$, and let $R = L^{-1}$ be the inverse of $L$. Lastly, let $M$ be the $N$ by $N$ tableau transfer matrix derived from $(Z, \omega, \pi, \theta)$, $S$, and $L$. Then for all integers $n \geq 1$ and for all $1 \leq i \leq N$ and $1 \leq j \leq N$, $M^n(i,j)$ equals to the number of elements $T$ of $$\Tb(\theta^{-n}(S) \cup \pi^{-1}([-n+1,0]), \omega)$$ that satisfy 
\begin{center}	
$T|_{S} \equiv_\theta R(i)$ \; and \; $T|_{\theta^{-n}(S)} \equiv_\theta R(j)$.
\end{center}
\end{corollary}

\begin{proof} This follows from the proof of Lemma \ref{tableau transfer matrices primitive}.
\end{proof}

A vector is \emph{non-negative} if all entries of the vector are non-negative and a vector is \emph{positive} if all entries of the vector are positive \cite{MA}. Moreover, let $M^T$ denote the transpose of a matrix $M$ and let $\langle u, v \rangle$ denote the standard inner product of vectors $u$ and $v$. Recall from spectral theory that the \emph{spectrum} $X \subset \mathbb{C}$ of a matrix $M$ is the set of distinct eigenvalues of $M$ and that the \emph{spectral circle} of a matrix $M$ is the set of distinct eigenvalues of $M$ of maximum modulus \cite{MA}. We state the part of the Perron-Frobenius theorem that we will need.

\begin{theorem}[Perron, Frobenius \cite{MA}]\label{Perron Frobenius} Let $M$ be a non-negative square matrix that is primitive. Then $M$ has exactly one eigenvalue $r$ on its spectral circle. Moreover, $r$ is a positive real number, $r$ is a simple eigenvalue of $M$, and an eigenvector $v$ of $M$ is non-negative if and only if $v$ is positive and $r$ is the eigenvalue corresponding to $v$.
\end{theorem}

Now, we prove the main result of this section.

\begin{theorem}\label{asymptotic formula} Let $(Z, \omega, \pi, \theta)$ be a periodic quadruple system. Moreover, let $P_n = \pi^{-1}([n])$ for all integers $n \geq 1$. Then there are constants $c > 0$ and $r > 1$ such that
$$|P_n, \omega| \sim c \, r^n$$
as $n \rightarrow \infty$.
\end{theorem}

\begin{proof} Let $(Z, \omega, \pi, \theta)$ be as described in Theorem \ref{asymptotic formula}. By Remark \ref{rem:separable index shapes always exist}, we can assume that $S$ is a separable index shape of $(Z, \omega, \pi, \theta)$. Next, let $L$ be an indexing of $\Tb(S, \omega)$. Moreover, let $n'$ be an admissible number for $(Z, \omega, \pi, \theta)$ and $S$, let $T_0 \in \Tb( \pi^{-1}([n',0]), \omega)$, and let $(v_n)_{n=1,2,\dots}$ be a sequence such that for all $n \geq 1$, $v_n$ is the $n^{th}$ vector derived from $(Z, \omega, \pi, \theta)$, $S$, $L$, and $T_0$. Lastly, let $M$ be the tableau transfer matrix derived from $(Z, \omega, \pi, \theta)$, $S$, and $L$. \\

Because $S$ is a separable index shape of $(Z, \omega, \pi, \theta)$, $M$ is separable so by Lemma \ref{tableau transfer matrices primitive}, $M$ is primitive. Hence, the Perron-Frobenius Theorem, and in particular Theorem \ref{Perron Frobenius}, applies to $M$. Let $r$ be the eigenvalue of $M$ in the spectral circle of $M$ and let $u$ be a positive eigenvector of $M$ that is an eigenvector of $r$. Lastly, let $N$ denote the number of rows of $M$. \\
	
Because $M$ is primitive, and because $r$ is as described in Theorem \ref{Perron Frobenius}, it follows from Theorem \ref{Perron Frobenius} and well-known results on spectral projectors of matrices (\cite{MA} p. 674, p. 630, p. 518) that there is a positive eigenvector $v$ of $M$ corresponding to $r$ and a positive eigenvector $u$ of $M^T$ corresponding to $r$ such that $\langle u, v \rangle > 0$ and
\begin{equation}\label{spectral projector}
\lim_{n \to \infty} M_0^n = \dfrac{v \, u^T}{\langle u, v \rangle} \,,
\end{equation}
where $M_0 = r^{-1} \, M$. \\

Consider the vector $v_1$. We show that $ \langle u, v_1 \rangle > 0$. By Lemma \ref{counting periodic shapes}, the sum of the entries of $v_1$ is $|P_1, \omega|$, where $P_1 = \pi^{-1}(\{1\})$. If $Q$ is a finite poset with at least one element and if $\omega_Q$ is a labelling of $Q$, then any \emph{injective} order-reversing map $f : Q \to \mathbb{N}_0$ is an element of $\mathcal{A}(Q, \omega_Q)$. In particular, $\mathcal{A}(P_1, \omega) \neq \emptyset$. So by Definition \ref{tableaux}, it follows that $|P_1, \omega| > 0$, implying that $v_1$ has at least one positive entry. Hence, as $u$ is positive, $\langle u, v_1 \rangle > 0$. \\

Write $M_1 = \lim_{n \to \infty} M_0^n$. Because $\langle u, v_1 \rangle > 0$, it follows from Equation \ref{spectral projector} that
$$ M_1 \, v_1 = \dfrac{v \, u^T \, v_1}{\langle u, v \rangle} = \dfrac{\langle u, v_1 \rangle}{ \langle u, v \rangle } \, v. $$
So as $\langle u, v_1 \rangle / \langle u, v \rangle > 0$, the following holds. Write $v_* = M_1 \, v_1$. Then by Theorem \ref{difference equation}, it follows from Equation \ref{spectral projector} that $v_n \sim r^n \, v_*$ as $n \to \infty$. Hence, the Theorem follows from Lemma \ref{counting periodic shapes}.
\end{proof}

A consequence of Theorem \ref{asymptotic formula} is the following, which adds to what is known about standard Young tableaux of parallelogramic shapes. Recall that $T(n,k)$ denotes the number of standard Young tableaux of parallelogramic shape $X_{n,k}$, where $X_{n,k}$ has $n$ rows and $k$ cells in each row.

\begin{corollary}\label{asymptotic SYT} Fix an integer $k \geq 3$. Then for some constants $c > 0$ and $r > 1$, $$T(n,k) \sim c \, r^n$$
as $n \to \infty$.
\end{corollary}

\begin{example} In \cite{CHC}, L\'opez et.al. computed the characteristic polynomials for the constant coefficient linear recurrence relations of minimum order for the sequences $(T(n,k))_{n=1,2,\dots}$ for $k$ up to twelve. From these polynomials, they computed closed-form expressions \cite{CHC} for $T(n,k)$ when $k = 4,5,6$. Using their closed-form expression for $k = 4$, we conclude that
$$T(n,4) \sim \bigg{(}\frac{1}{2} - \frac{1}{2 \, \sqrt{2}}\bigg{)} \, (3 + 2 \sqrt{2})^n$$
as $n \to \infty$.
\end{example}

\begin{proof} As described in the proof of Corollary \ref{existence}, Proposition \ref{Full Euclidean system} implies that collections of parallelogramic shapes constitute examples of periodic quadruple systems if the number of cells in each row is fixed. Hence, the corollary follows from Theorem \ref{asymptotic formula}.
\end{proof}

Another corollary of Theorem \ref{asymptotic formula} is the following. Recall that $X_{n,k}$ denotes a parallelogramic shape with $n$ rows and $k$ cells in each row.

\begin{corollary} For integers $n \geq 1$ and $k \geq 3$, let $T'(n,k)$ denote the number of semi-standard tableau classes of shape $X_{n,k}$. Then for some constants $c > 0$ and $r > 1$, $$T'(n,k) \sim c \, r^n$$
as $n \to \infty$.
\end{corollary}

\begin{proof} The proof is similar to the proof of Corollary \ref{asymptotic SYT}.
\end{proof}

Lastly, we remark on another application of the results in this section. Let $\tr M$ denote the trace of a square matrix $M$.

\begin{remark} If we enumerate the sequences in Theorem \ref{recurrence relations} with \emph{characteristic polynomials} of tableau transfer matrices instead of minimal polynomials of transfer transfer matrices, then we can describe the coefficients of the resulting constant coefficient linear recurrence relations. Recall that Corollary \ref{powers} gives combinatorial descriptions for the entries of arbitrary powers of transfer matrices. In particular, using the notation in Corollary \ref{powers}, the traces of these matrices, for all integers $n \geq 1$, can be combinatorially described as being the number of elements $T$ in 
$$\Tb(\theta^{-n}(S) \cup \pi^{-1}([-n+1, 0]), \omega)$$
such that $T|_{\theta^{-n}(S)} \equiv_\theta T|_S$. Hence, we can describe the aforementioned coefficients using known results for describing the coefficients of a characteristic polynomial. For instance, we can use Lewin's result \cite{Lewin} that expresses the coefficients of the characteristic polynomial of a square matrix $M$ solely in terms of parameters of the form $\tr M^n$ and that involves a summation over the (integer) partitions of a fixed integer.
\end{remark}


\section*{Acknowledgements}\label{sec:acknow} The author would like to thank Stephanie van~Willigenburg for her guidance and advice during the development of this paper.
\bibliographystyle{amsplain}

\end{document}